\documentclass[11pt]{article}

\usepackage{amsmath}
\usepackage{graphics}
\usepackage{enumerate}
\usepackage{epsfig}
\usepackage{amsfonts}
\usepackage{amssymb}

\begin{document}

\def\ALERT#1{{\large\bf $\clubsuit$#1$\clubsuit$}}

\newtheorem{defin}{Definition}
\newtheorem{theorem}{Theorem}
\newtheorem{notice}{Notice}
\newtheorem{lemma}{Lemma}
\newtheorem{cor}{Corollary}
\newtheorem{example}{Example}
\newtheorem{remark}{Remark}
\def\begproof{\noindent{\bf Proof: }}
\def\endproof{\quad\vrule height4pt width4pt depth0pt \medskip}
\def\div{\nabla\cdot}
\def\rot{\nabla\times}
\def\sign{{\rm sign}}
\def\arsinh{{\rm arsinh}}
\def\arcosh{{\rm arcosh}}
\def\diag{{\rm diag}}
\def\const{{\rm const}}
\def\d{\,\mathrm{d}}
\def\rho{\varrho}
\def\eps{\varepsilon}
\def\phi{\varphi}
\def\theta{\vartheta}
\def\N{\mathbb{N}}
\def\R{\mathbb{R}}
\def\C{\hbox{\rlap{\kern.24em\raise.1ex\hbox
      {\vrule height1.3ex width.9pt}}C}}
\def\Q{\hbox{\rlap{\kern.24em\raise.1ex\hbox
      {\vrule height1.3ex width.9pt}}Q}}
\def\M{\hbox{\rlap{I}\kern.16em\rlap{I}M}}
\def\Z{\hbox{\rlap{Z}\kern.20em Z}}
\def\({\begin{eqnarray}}
\def\){\end{eqnarray}}
\def\[{\begin{eqnarray*}}
\def\]{\end{eqnarray*}}
\def\part#1#2{{\partial #1\over\partial #2}} 
\def\partk#1#2#3{{\partial^#3 #1\over\partial #2^#3}} 
\def\mat#1{{D #1\over Dt}}
\def\dx{\nabla_x}
\def\dv{\nabla_v}
\def\grad{\nabla}
\def\tT{{\mbox{\tiny{T}}}}

\def\Norm#1{\left\| #1 \right\|}
\def\pmb#1{\setbox0=\hbox{$#1$}
  \kern-.025em\copy0\kern-\wd0
  \kern-.05em\copy0\kern-\wd0
  \kern-.025em\raise.0433em\box0 }
\def\bar{\overline}
\def\lbar{\underline}
\def\fref#1{(\ref{#1})}
\def\half{\frac{1}{2}}
\def\oo#1{\frac{1}{#1}}

\def\tot#1#2{\frac{\d #1}{\d #2}} 
\def\laplace{\Delta}
\def\d{\,\mathrm{d}}
\def\N{\mathbb{N}}
\def\R{\mathbb{R}}
\def\supp{\mbox{supp }}
\def\rho{\varrho}
\def\eps{\varepsilon}
\def\phi{\varphi}

\def\E{\mathcal{E}}
\def\K{\mathcal{K}}
\newcommand{\zb}[1]{\ensuremath{\boldsymbol{#1}}}

\title{Consistency of Variational Continuous-Domain Quantization\\
 via Kinetic Theory}
\author{
Massimo Fornasier\footnote{Faculty of Mathematics, Technical University of Munich, Boltzmannstrasse 3, D-85748 Garching, Germany,
email: {\tt massimo.fornasier@ma.tum.de}.}, \
Jan Ha\v{s}kovec\footnote{Johann Radon Institute for Computational and
Applied Mathematics, Austrian Academy of Sciences, Altenbergerstrasse 69, A-4040 Linz, Austria,
email: {\tt  jan.haskovec@oeaw.ac.at}.} \ and
Gabriele Steidl\footnote{Department of Mathematics, Technical University of Kaiserslautern, Germany
email: {\tt  steidl@mathematik.uni-kl.de}.}
}
\maketitle

\begin{abstract}
We study the kinetic mean-field limits of the discrete systems of interacting particles
used for halftoning of images in the sense of continuous-domain quantization. 
Under mild assumptions on the regularity of the interacting kernels we provide a rigorous derivation
of the mean-field kinetic equation. Moreover, we study the energy of the system,
show that it is a Lyapunov functional and prove that in the long time limit
the solution tends to an equilibrium given by a local minimum of the energy.
In a special case we prove that the equilibrium is unique and is identical
to the prescribed image profile. This proves the consistency of the particle halftoning method
when the number of particles tends to infinity.
\end{abstract}

\noindent
{\bf AMS subject classification (MSC 2010):} 82C10, 82C22, 68U10
\\



\section{Introduction} \label{sec:introduction}
\begin{figure}[ht]
{\centering \begin{tabular}[h]{lr}
\resizebox*{0.40\linewidth}{!}{\includegraphics{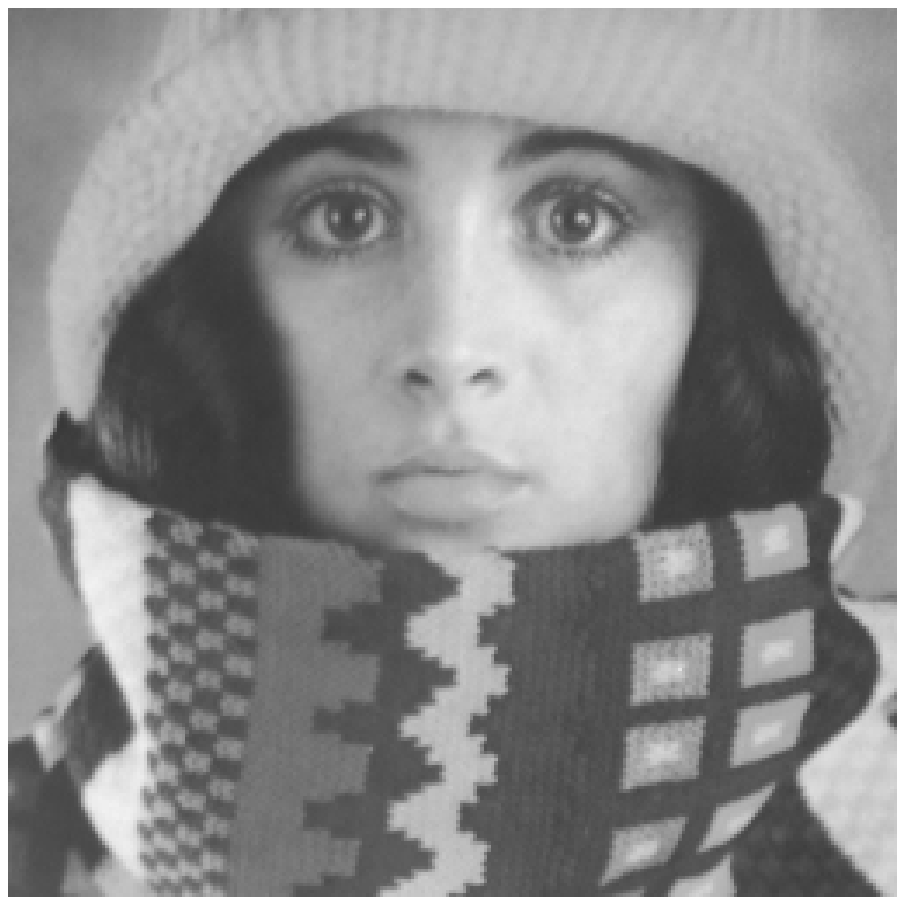}} $\quad$ &
$\quad$ \resizebox*{0.40\linewidth}{!}{\includegraphics{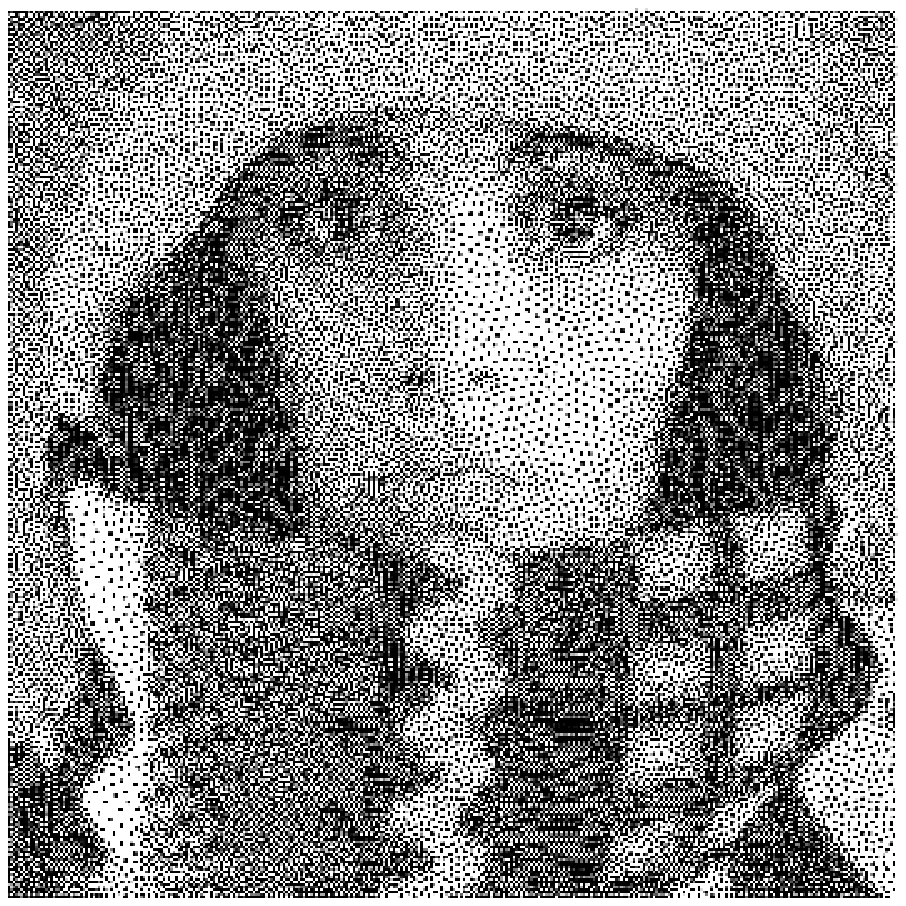}}
\end{tabular}\par}
\caption{\label{fig:1}
Left: Original $256 \times 256$ image.
Right: Halftoning result with $m = 30150$ points using the technique from \cite{TSGSW}.
}
\end{figure}

A halftoning  method places black dots in an image
in such a way that their density gives the impression of tone.
For an illustration see Fig. \ref{fig:1}.
Due to its various applications, halftoning is an active field of research
and we refer to the recent papers \cite{SGBW10, Fa11} for deterministic and, resp., stochastic point distributions. 

In this paper, we consider a continuous-domain quantization method based on electrostatic-like principles 
studied in~\cite{TSGSW}, where the basic idea goes back to \cite{SGBW10}. 
The method consists 
in considering a system of $N$ particles with electrostatic-like interaction (repulsion)
and exposed to an attractive external potential $w$ with compact support in $\R^d$  which represents the image 
to be approximated by points.
In \cite{TSGSW} the authors considered the discrete energy functional
\begin{equation} \label{starting point}
    E(p) := \sum_{k=1}^N \int_{\R^d} w(x) |p_k-x| \d x - \frac{\lambda}{2} \sum_{k=1}^N \sum_{l=1}^N |p_k-p_l |, 
\end{equation}
where $p:=(p_1, \dots, p_N)^\tT \in \mathbb R^{N\times d}$, the Euclidean distance is denoted  by $| \cdot|$,
\begin{equation} \label{eq0}
   \lambda := \frac{1}{N} \int_{\R^d} w(x) \d x \,,
\end{equation}
and defined the evolution of the particle system as the corresponding gradient flow
$\partial_t p \in - \partial_p E$.
They showed that for $d=1$ the energy functional is continuous and coercive,
and calculated explicitly its minimizers. In the two-dimensional setting
it is not possible to obtain explicit expressions for the minimizers anymore.
Instead, the authors employed a difference of convex functions (DC) algorithm
together with fast summation methods for non-equispaced knots
to obtain a local minimum of the variational problem numerically.
In \cite{GPS11} the function $w$ was also considered on other sets such as 
$\mathbb T^d$ or $\mathbb S^2$ 
and kernels other than the Euclidian distance were used. 
This generalized approach stems from the study of optimal quadrature error functionals 
on reproducing kernel Hilbert spaces with respect to the quadrature knots and is closely related to
so-called {\it discrepancy functionals}~\cite{Graf-Luschgy}.

In this paper we study the kinetic mean-field limits of discrete systems like \eqref{starting point},
that are obtained as the number of particles tends to infinity.
First, in Section~\ref{sec:kernel-based}, we specify the interaction kernels which we will consider
and provide some particular examples. Moreover, we introduce a generalized setting
with different kernels for the attractive and repulsive interactions,
and show that this new setting can be reduced to the previous one with appropriately modified data.
In Section~\ref{sec:mean-field-limit} we provide, under mild regularity conditions on the interaction kernels,
a rigorous derivation of the mean field kinetic equation obtained in the limit as the number of particles tends to infinity.
We show that the corresponding energy functional, obtained as a formal limit
of the discrete energy, is a Lyapunov functional for the kinetic equation and that in the long-time
limit the solution tends to an equilibrium, which is a local minimum of the energy.
For a special choice of the interaction kernel, we are able to show that the equilibrium
coincides with the prescribed image profile. This proves the consistency of the discrete
halftoning method when the number of particles tends to infinity.
Finally, in Section~\ref{sec:numerics}, we provide numerical examples confirming our theoretical results, and showing the behavior
of the model subject to different competitive attraction and repulsion terms as well as the
consistency of the particle system with respect to its continuous kinetic limit.

%
\section{Energy Functionals with Kernels} \label{sec:kernel-based}
In the following, let $\Omega$ be a domain, 
which we will be either the Euclidean space $\Omega=\R^d$
or the torus $\Omega=\mathbb T^d$,  $d\geq 1$.
A symmetric function $K: \Omega \times \Omega \rightarrow \R$ is said to be \emph{positive semi-definite}
if for any $N \in \N$ points $x_1, \dots, x_N \in \Omega$ and any $a\in\R^N\setminus\{0\}$ the
relationship $a^T (K(x_i, x_j ))_{i,j=1}^N a \geq 0$ holds true,
and \emph{positive definite} if we have strict inequality.
Let $K: \Omega \times \Omega \rightarrow \R$ be a symmetric, positive semi-definite function, and $H_K$ be the reproducing kernel
Hilbert space (RKHS) associated with $K$, see \cite{Ar50}.
Then we are interested in the functional
\begin{equation} \label{thats_it}
E_K (p) := \sum_{i=1}^N \int_{\Omega} w( x) K( p_i, x) \, \d  x
  - \frac{\lambda}{2} \sum_{i,j=1}^N  K( p_i, p_j) \,.
\end{equation}
In the case $\Omega = \R^d$ we suppose that the function $w:\R^d \rightarrow \R_+$ is compactly supported.
In \cite{GPS11} it was shown that this functional is related to the optimality
of a certain quadrature rule for functions in $H_K$
depending on the knots $p_i$, $i=1,\ldots,N$. 
By the following remark, which was proved in \cite{GPS11}, slight modifications of the kernel do not change the minimizers of the functional $E_K$.
%
\begin{remark} \label{rem:1}
Let $K: \Omega \times \Omega \rightarrow \mathbb R$ be a symmetric function and
$\tilde K(x,y) := a K(x,y) + b(K(x,0) + K(0,y)) + c$ 
with $a > 0$ and $b,c \in \mathbb R$.
Then the minimizers of $E_K$ and $E_{\tilde K}$ coincide.
\end{remark}
%
In this paper, we consider
{\it radial kernels} on $\Omega = \mathbb \R^d$, i.e., 
\begin{equation} \label{rad_kernel}
K(x,y) = \varphi(x-y) = \Phi(|x-y|)
\end{equation}
with $\Phi:[0,\infty) \rightarrow \mathbb R$.
For the 1-periodic setting $\Omega = \mathbb T^1$ we use the same notation,
where $| x-y | $ has to be replaced by the ``periodic'' distance $\min\{|x-y|, 1-|x-y|\}$. 
In the case $\Omega = \mathbb T^d$ we consider tensor products of radial kernels.
We call $\varphi$ positive semi-definite (resp. positive definite) 
if the corresponding kernel is positive semi-definite (resp. positive definite).
The functional of interest becomes then
\begin{equation} \label{kernel_phi}
   E_\varphi (p) := \sum_{i=1}^N \int_{\Omega} w(x) \phi(p_i - x) \, \d x
       - \frac{\lambda}{2} \sum_{i,j=1}^N  \varphi(p_i - p_j) \,.
\end{equation}
%
{\bf Example.}
We give some interesting examples of positive semi-definite kernels, see \cite{We05} and \cite{Wa90}.
\begin{enumerate}
\item Let ${\Omega} := \mathbb R^d$. 
Then the functions
\begin{eqnarray*}
 \varphi(x) &=& (1-|x|)^\tau_+, \quad \tau \ge \left\lfloor \frac{d}{2} \right\rfloor + 1,\\
 \varphi(x) &=& (\varepsilon^2 + |x|^2)^{-\beta}, \quad \beta > \frac{d}{2} \qquad{\rm (inverse \; multiquadrics)}
\end{eqnarray*}
are positive definite.
Next, consider the {\it conditionally positive definite radial kernels of order 1} defined by
\begin{eqnarray*}
\varphi(x) &:=& - |x|^\tau,
\quad 0< \tau < 2, \\
\varphi(x) &:=& - (\varepsilon^2 + |x|^2)^\tau,
\quad 0 < \tau< 1, \quad
{\rm (multiquadrics)}. 
\end{eqnarray*}
The kernels $K(x,y) = \varphi(x-y)$ are not positive semi-definite.
However, their slight modifications  given by
\begin{equation*} 
\tilde K(x,y) := \varphi(x-y)  - \varphi(y) - \varphi(x) + \varphi(0) 
\end{equation*}
define positive semi-definite kernels, and the corresponding  RKHSs were characterized in \cite[Theorem 10.18 ]{We05}.
By Remark \ref{rem:1}, $E_K$ and $E_{\tilde K}$
have the same minimizers, so that we can work with the original kernel $K$ in the energy functional. 
\item
Let $\Omega := \mathbb T^1$.
Up to an additive constant, Wahba's spline kernels are given by
$$
K(x,y) := \frac{(-1)^{m-1}}{(2m)!} B_{2m} (|x-y|) = 2 \sum _{k=1}^\infty \frac{1}{(2 \pi k)^{2m}} \cos(2 \pi k (x-y))
$$
where $B_{2m}$ denotes the {\it Bernoulli polynomial} of degree $2m$. 
Note that $B_{2m} (1-t) = B_{2m} (t)$.
For example, we have
$$
B_2(t) = t^2 - t + \frac16,\qquad
B_4(t) = t^4 - 2t^3 + t^2-\frac{1}{30}.
$$
\end{enumerate}
%
Let us now introduce a slight generalization of~\eqref{kernel_phi},
where we consider different kernels for the attractive interaction
(attraction of the particles by the image profile $w$)
and the repulsive interaction (particle-particle repulsion).
In particular, we introduce the functions $\varphi$ and $\psi$ 
related to different radial kernels \eqref{rad_kernel},
and the generalized energy functional
\begin{equation} \label{kernel_psi_phi}
   E_{\psi,\varphi} (p) := \sum_{i=1}^N \int_{\Omega} w(x) \phi(p_i - x) \, \d x
         - \frac{\lambda}{2} \sum_{i,j=1}^N  \psi(p_i - p_j) \,.
\end{equation}
The following remark gives an intuition on the behavior of the corresponding quantization process.

\begin{remark} \label{rem:2}
Assume that the kernel $K$ in \eqref{thats_it}
is in addition continuous and an element of $L^2({\Omega} \times {\Omega})$ (Mercer kernel). 
Then it can be expanded into an absolutely and uniformly convergent series,
\begin{equation*} \label{mercer_kernel}
K(x,y) = K_\lambda (x,y):= \sum_{\ell=1}^\infty \lambda_\ell \eta_\ell(x) \eta_\ell(y) 
\end{equation*}
of orthonormal eigenfunctions $\eta_\ell \in L^2(\Omega)$ and associated eigenvalues $\lambda_\ell >0$
of the compact, self-adjoint integral operator $T_K$ on $L^2(\Omega)$ given by
$$
T_K f(x) := \int_{\Omega} K(x,y) f(y) \, \d  y.
$$
Assume further that $w$ can be also expanded into an absolutely convergent series
$
w(x) := \sum_{k=1}^\infty w_k \eta_k(x).
$
Then the functional \eqref{thats_it} becomes
\begin{eqnarray*}
E_{K} (p) 
&=&
\sum_{i=1}^N \int_{\Omega} w( x) K_\lambda( p_i, x) \, \d  x - \frac{\lambda}{2} \sum_{i,j=1}^N  K_\lambda( p_i, p_j)\\
&=&  
\sum_{i=1}^N \int_{\Omega} \sum_{k=1}^\infty w_k \eta_k(x) \, \sum_{\ell=1}^\infty \lambda_\ell \eta_\ell(p_i) \eta_\ell(x) \d x
-\frac{\lambda}{2} \sum_{i,j=1}^N  K_\lambda ( p_i, p_j) \\
&=&
\sum_{i=1}^N \sum_{k,\ell=1}^\infty  w_k \lambda_\ell \eta_\ell(p_i) \int_{\Omega} \eta_k(x) \eta_\ell(x) \d x
- \frac{\lambda}{2} \sum_{i,j=1}^N  K_\lambda ( p_i, p_j)\\
&=&
\sum_{i=1}^N  \sum_{k=1}^\infty w_k \lambda_k \eta_k(p_i)
- \frac{\lambda}{2} \sum_{i,j=1}^N  K_\lambda( p_i, p_j) \,.
\end{eqnarray*}
On the other hand, if we consider $E_K$ for another function $\tilde w(x) := \sum_{k=1}^\infty w_k v_k \eta_k(x)$,
we have
\begin{eqnarray*}
E_K (p) 
&=&
\sum_{i=1}^N \int_\Omega \tilde w(x) K_\lambda (p_i,x) \, \d x
- \frac{\lambda}{2} \sum_{i,j=1}^N  K_\lambda( p_i, p_j) \\
&=& 
\sum_{i=1}^N  \sum_{k=1}^\infty w_k v_k \lambda_k \eta_k(p_i)
- \frac{\lambda}{2} \sum_{i,j=1}^N  K_\lambda( p_i, p_j)\\
&=&
\sum_{i=1}^N \int_{\Omega} w(x) K_\mu( p_i, x) \, \d  x
- \frac{\lambda}{2} \sum_{i,j=1}^N  K_\lambda( p_i, p_j)
\end{eqnarray*}
with $K_\mu(x,y) := \sum_{\ell=1}^\infty \mu_\ell \eta_\ell(x) \eta_\ell(y)$
and $\mu_k := v_k\lambda_k$,
where we assume absolute convergence of the involved series. Hence,
using a smoother kernel $K_\mu$ for the interaction with the datum $w$
than $K_\lambda$ (i.e., $\mu_k$ decays faster than $\lambda_k$)
leads to the approximation of a smoother function $\tilde w$
($w_k v_k$ decays faster than $w_k$), and vice versa.
\end{remark}

\section{Mean-Field Limit} \label{sec:mean-field-limit}

We are interested in the passage to the limit when the number of particles $N$ tends to infinity.
For simplicity, we restrict our attention to radial kernels  and $\Omega = \mathbb R^d$
although the analysis works as well for the periodic settings  $\Omega = \mathbb T$
and $\Omega = \mathbb T^d$, $d \ge 2$ with the tensor product of radial kernels.
Moreover, without loss of generality, we prescribe the normalization 
$$\int_{\Omega} w(x) \d x =1, \quad N \lambda  = 1$$
and suppose that $w \ge 0$ is compactly supported.

\subsection{Passage to the Mean-Field Limit} 
The evolution of the $N$-particle system according
to the gradient flow of the discrete generalized energy functional 
\( \label{discrete energy}
   E(p) = E_{\varphi,\psi} (p) := \sum_{k=1}^N \int_\Omega w(x) \phi(p_k-x) \d x - \frac{1}{2N} \sum_{k,\ell=1}^N  \psi(p_k-p_\ell)
\)
is given, under the assumption that $\phi, \psi\in C^1(\Omega)$, by
\begin{eqnarray}
    \tot{}{t} p_i(t) &=& - \grad_{p_i} E(p(t)) \nonumber\\
&=& - \int_\Omega w(x) \grad\phi(p_i(t)-x) \d x + \frac{1}{N} \sum_{\ell=1}^N \grad \psi(p_i(t)-p_\ell(t))\,,
     \qquad i=1,\dots,N , \label{ODE}
\end{eqnarray}
subject to the initial condition
\(
   p_i(0) = p_i^0 \,,\qquad i=1,\dots,N. \label{ODE_IC}
\)
The mean field limit is obtained as the number of particles $N$
tends to infinity. Then, the vector of time-dependent particle
positions $p(t)\in\Omega^N$ is replaced by the time-dependent probability measure
$f(x,t)$, where, roughly speaking, $f(x,t) \d x$
can be understood as the probability that a particle is located
in the space element $\d x$ around the position $x\in\Omega$
at time $t\geq 0$.

In the following, let $\mathcal M(\Omega)$ 
denote the space of Radon measures on $\Omega$
and $C_c(\Omega)$ the space of continuous, compactly supported functions on $\Omega$.
Further, let
$L^\infty(\mathbb R_+,\mathcal M(\Omega))$ 
denote the space of functions from $\mathbb R_+$ to $\mathcal M(\Omega)$ which are essentially bounded,
i.e., $f: t \rightarrow f(\cdot,t) = f_t$
with ${\rm ess} \, {\rm sup}_{t \in \mathbb R_+}  \int_{\Omega} \d |f(\cdot,t)| < \infty$.
Note that $L^\infty(\mathbb R_+,\mathcal M(\Omega))$ is the dual space of
$L^1 (\mathbb R_+,C_c(\Omega))$, 
the space of functions from $\mathbb R_+$ to $C_c(\Omega)$
such that
$\int_0^\infty  \|g(\cdot,t) \|_\infty \d t < \infty$, see, e.g.,~\cite{Bourbaki}.
Moreover, let us denote by $\mathcal M^1(\Omega)$ the set
of probability measures on $\Omega$, i.e., $f\in\mathcal M^1(\Omega)$
if and only if $f$ is a nonnegative Radon measure such that $\int_\Omega \d f = 1$.

For any $N\in\N$, let us denote by $f^N$ the empirical measures
\(   \label{empirical}
   f^N(\cdot,t) = \frac{1}{N} \sum_{i=1}^N \delta(\cdot - p_i(t))
\)
corresponding to the evolution of the $N$-particle system \eqref{ODE}.
Then, each $f^N$ is a time-dependent probability measure,
such that $f^N \in L^\infty(\mathbb R_+,\mathcal M(\Omega))$.
In the following theorem we carry out the rigorous mean field limit passage $N\to\infty$.
%
\begin{theorem} \label{theorem:mean-field}
Let $\varphi, \psi \in C^1(\Omega)$, where $\grad\psi$ is in addition bounded.
Let  $(f^N)_{N\in\N}$ be given by~\eqref{empirical}, corresponding to the system \eqref{ODE}
with the initial datum~\eqref{ODE_IC}.
Moreover, assume that there exists a probability measure $f_0\in\mathcal M^1(\Omega)$ such that
$f^N(\cdot,0) \to f_0(\cdot)$ weakly-* in $\mathcal M(\Omega)$ as $N\to\infty$.

Then there exists a subsequence $\left( f^{N_k} \right)_{k\in\N}$
which converges weakly-* in $L^\infty(\mathbb R_+,\mathcal M(\Omega))$ to
a time-dependent probability measure
$f \in L^\infty(\mathbb R_+,\mathcal M^1(\Omega))$ which solves, in the sense of distributions,
the mean-field equation
\begin{eqnarray}  
    \partial_t f 
     &=& \grad_y \cdot\left( \int_{\Omega} (w(x)\grad\phi(y-x)
               - f(x,t)\grad\psi(y-x)) f(y,t) \d x \right) \nonumber \\
     &=& \nabla \cdot \left( \grad\K[f] f \right) \,, \label{kinetic}                
\end{eqnarray}
where
\( \label{K}
   \K[f](y,t) := \int_{\Omega} \left( w(x) \phi(y-x) - f(x,t) \psi(y-x) \right) \d x \,,
\)
subject to the initial condition
\(   \label{kinetic_IC}
    f(\cdot,0) = f_0 \,.
\)
\end{theorem}
\begproof
First we show that for all $N\in\N$ the empirical measures~\eqref{empirical}
are distributional solutions of~\eqref{kinetic}.
Indeed, considering a smooth, compactly supported test function $\xi \in C^\infty_c(\Omega\times[0,\infty))$,
we obtain
\[
  && \int_0^\infty \int_{\Omega} f^N(y,t) \partial_t \xi(y,t) \d y \d t + \int_{\Omega} f^N(y,0) \xi(y,0) \d y \\
  &=& \frac{1}{N} \int_0^\infty\sum_{i=1}^N \partial_t \xi(p_i(t),t) \d t + \frac1N \sum_{i=1}^N \xi(p_i(0),0) \\
  &=& \frac{1}{N} \int_0^\infty\sum_{i=1}^N \left[\tot{}{t} \xi(p_i(t),t) - \grad_{p_i}\xi(p_i(t),t)\cdot\tot{}{t}p_i(t) \right] \d t + \frac1N \sum_{i=1}^N \xi(p_i(0),0) \\
  &=& \frac{1}{N} \int_0^\infty\sum_{i=1}^N \grad_{p_i} \xi(p_i(t),t) \cdot \left[ \int_{\Omega} w(x) \grad\phi(p_i(t)-x) \d x- \frac{1}{N} \sum_{\ell=1}^N \grad \psi(p_i(t)-p_\ell(t)) \right] \d t \\
  &=& \int_0^\infty\int_{\Omega} f^N(y,t) \grad_y\xi(y,t) \cdot \left[ \int_{\Omega} w(x) \grad\phi(y-x) \d x - 
	\int_{\Omega} f^N(x,t) \grad \psi(y-x) \d x \right] \d y \d t \,,
\]
Therefore, we have the identity
\(
  && \int_0^\infty \int_{\Omega} f^N(y,t) \partial_t \xi(y,t) \d y \d t + \int_{\Omega} f^N(y,0) \xi(y,0) \d y \label{zwischen} \\
  &=& \int_0^\infty\int_{\Omega} f^N(y,t) \grad_y\xi(y,t) \cdot \left[ \int_{\Omega} w(x) \grad\phi(y-x) \d x - 
	\int_{\Omega} f^N(x,t) \grad \psi(y-x) \d x \right] \d y \d t \,, \nonumber
\)
for all test functions $\xi \in C^\infty_c(\Omega\times[0,\infty))$,
that is the distributional formulation of~\eqref{kinetic} with $f^N$ in place of $f$
and the initial condition $f^N(\cdot,0) = \frac{1}{N} \sum_{i=1}^N \delta(\cdot - p_i^0)$.
                                              
Now, since $(f^N)_{N\in\N}$ is a sequence of time-dependent probability measures,
it is uniformly bounded in $L^\infty(\mathbb R_+,\mathcal M(\Omega))$,
so that there exists a subsequence $(f^{N_k})_{k\in\N}$ which converges weakly-* to some $f$
in $L^\infty(\mathbb R_+,\mathcal M^1(\Omega))$.
We show that $f$ is a distributional solution of~\eqref{kinetic}.
The limit passage in the linear terms of \eqref{zwischen} follows immediately.
Moreover, since $\grad\psi$ is assumed to be continuous and bounded, the sequence
$\int_{\Omega} f^N(x,t) \grad \psi(y-x) \d x$ is uniformly equi-continuous and uniformly bounded
for $y\in\supp\xi(\cdot,t)$ and for almost every $t\in\R_+$.
Therefore, due to the Arzel\`a-Ascoli theorem, this sequence converges strongly in $L^\infty(\supp\xi(\cdot,t))$
for almost all $t\in\R_+$. Finally, the bounded convergence theorem ensures the strong convergence
of the sequence in $L^1(\R_+, L^\infty(\supp\xi))$ and this justifies the limit passage in the nonlinear term.

Therefore, in the limit $N\to\infty$, we have obtained
\(
  && \int_0^\infty \int_{\Omega} f(y,t) \partial_t \xi(y,t) \d y \d t + \int_{\Omega} f_0(y) \xi(y,0) \d y \label{weak} \\
  &=& \int_0^\infty\int_{\Omega} f(y,t) \grad_y\xi(y,t) \cdot \left[ \int_{\Omega} w(x) \grad\phi(y-x) \d x - 
	\int_{\Omega} f(x,t) \grad \psi(y-x) \d x \right] \d y \d t \,, \nonumber
\)
which is the distributional formulation of~\eqref{kinetic} subject to the initial condition~\eqref{kinetic_IC}.
\endproof

Equation \eqref{kinetic} describes the evolution of the time-dependent probability measure $f(\cdot,t)$
due to the mutual repulsive interaction between the particles and the
attractive interaction with the datum $w$.
In the following lemma we show that, under mild regularity assumptions on $\phi$, $\psi$ and $f_0$,
the solution $f$ of \eqref{kinetic} is in fact classical.

\begin{lemma} \label{lemma_1}
Let $\grad\phi$ and $\grad\psi$ be globally Lipschitz continuous on $\mathbb R^d$, 
i.e., there exist constants $L_1,L_2$ such that
\[
    |\grad\phi(x) - \grad\phi(y)| &\leq& L_1 |x-y| \,,\\
    |\grad\psi(x) - \grad\psi(y)| &\leq& L_2 |x-y| \,,
\]
for all $x,y \in \Omega$.
Let $f_0 \in C_c^1(\Omega)$ be nonnegative, compactly supported and fulfill $\int_\Omega f_0(x) \d x = 1$.
Then the corresponding distributional solution $f$ of~\eqref{kinetic}
subject to the initial condition \eqref{kinetic_IC} is in fact a classical solution 
with $f \in C^1(\Omega\times\mathbb R_+)$
and $f(\cdot,t) \geq 0$ for all $t\geq 0$.
Moreover, $f(\cdot,t)$ is compactly supported on $\Omega$ for any $t \in \mathbb R_+$.
\end{lemma}
\begproof
Since the distributional solution
$f \in L^\infty(\mathbb R_+,\mathcal M^1(\Omega))$ constructed in Theorem~\ref{theorem:mean-field}
is a time-dependent probability measure, we have
$f(\cdot,t) \ge 0$ and
\[
   \int_{\Omega} \d f(t,\cdot) \equiv 1\quad\mbox{for all } t\geq 0 \,.
\]
The essential point is to observe that due to the assumptions on $\phi$ and $\psi$,
the transport field $\grad\K[f]$ for~\eqref{kinetic}
is Lipschitz continuous. Indeed, we have
\[
    |\grad\K[f](p) - \grad\K[f](q)| &\leq&
      \int_{\Omega} w(x) |\grad\phi(p-x) - \grad\phi(q-x)| \d x \\
      &&+
      \int_\Omega f(x,t) |\grad\psi(p-x) - \grad\psi(q-x)| \d x \\
    &\leq&
      L_1\int_\Omega w(x) |p-q| \d x +
       L_2\int_\Omega f(x,t) |p-q| \d x  \\
    &\leq&
      (L_1+L_2) |p-q| \,.
\]
Therefore, $f$ is a solution of a hyperbolic transport equation
with globally Lipschitz-continuous transport field $\grad\K[f]$.
As such, the values of the initial condition $f_0$ propagate along the characteristics
\(    \label{characteristics}
    \dot x(t) = \grad\K[f](x(t))
\)
with finite speeds.
Due to the assumption $f_0 \in C_c^1(\Omega)$,
from the standard theory of hyperbolic transport equations
(method of characteristics and Cauchy-Lipschitz theorem, see~\cite{Evans})
it follows that $f\in C^1(\Omega\times\mathbb R_+)$.
The compactness of the support of $f$ for all times follows from
the assumed compactness of the support of $f_0$ and
the finite characteristic speeds~\eqref{characteristics}.
\endproof

\begin{remark}\label{rem:3}
Similarly as in Remark~\ref{rem:2}, we observe that the evolution induced by~\eqref{kinetic}
with two different interaction kernels $\phi$ and $\psi$ is in fact equivalent
to an evolution produced by using same interaction kernels, but with a modified data $w$.
Indeed, assuming that $\phi$ and $\psi$ are continuously differentiable
and that the Fourier transform of $\psi$ is nonzero almost everywhere
(as is the case for positive definite kernels),
the Fourier transform of~\eqref{kinetic} reads
\[
    \partial_t \hat f &=& i\xi\cdot \left( (\hat w\widehat{\grad\phi} - \hat f\widehat{\grad\psi})\ast\hat f\right) \\
              &=& - |\xi|^2 \left( (\hat w\hat{\phi} - \hat f\hat{\psi})\ast\hat f\right) \\
              &=& - |\xi|^2 \left( \hat{\psi} \left(\hat w\frac{\hat{\phi}}{\hat\psi} - \hat f\right)\ast\hat f\right) \,.
\]
Applying the inverse Fourier transform, we get
\[
    \partial_t f = \grad\cdot\left((\tilde w-f)\ast\grad\psi)f\right) \,,
\]
where $\tilde w$ is the inverse Fourier transform of $(\hat w\hat{\phi}/\hat\psi)$.
Therefore, taking $\phi$ smoother than $\psi$ corresponds to a smoothing of $w$,
while $\psi$ smoother than $\phi$ corresponds to a ``sharpening'' (anti-smoothing) of~$w$.
\end{remark}

\subsection{Energy dissipation, long time behavior and equilibria}
Let us observe that the formal limit
of the discrete energy~\eqref{discrete energy} as $N\to\infty$
is given by the continuous energy functional
\(  \label{cont energy}
    \E[f] =  \int_\Omega \int_\Omega w(x) \phi(p-x) f(p) \d x\d p
            - \frac{1}{2} \int_\Omega\int_\Omega f(x) \psi(p-x) f(p) \d p\d x  \,,
\)
defined for all $f\in\mathcal M^1(\Omega)$.
Let us mention that the corresponding formal gradient flow with respect to the topology
of 2-Wasserstein distance on the space of probability measures (see~\cite{JKO} or~\cite{Otto} for details)
is given by the hyperbolic transport equation \eqref{kinetic}.
This suggests that the energy $\E[f(x,\cdot)]$
actually is a Lyapunov functional, thus nonincreasing along the solutions of~\eqref{kinetic}:
\(    \label{energy dissipation}
    \tot{}{t} \E[f(x,t)] \leq 0 \,.
\)
Indeed, at least for classical solutions, this inequality can be proven rigorously:

\begin{lemma} \label{lemma:classical_soln}
Let $\varphi, \psi$ be of the form~\eqref{rad_kernel}.
Let $f\in C^1(\Omega\times\mathbb R_+)$ be a classical solution of~\eqref{kinetic}--\eqref{kinetic_IC}
in the sense of Lemma~\ref{lemma_1}.
Then~\eqref{energy dissipation} holds true.
\end{lemma}

\begproof
The proof follows from the direct calculation
\begin{eqnarray}
   \tot{}{t} \E[f(x,t)]
   &=& \int_\Omega\int_\Omega w(x) \varphi(p-x) \nabla_p \cdot \left( \grad\K[f](p,t) f(p,t) \right)  \d p\d x \nonumber \\
   && - \frac12 \int_\Omega\int_\Omega \psi(p-x) f(x,t) \nabla_p \cdot \left(\grad\K[f](p,t) f(p,t) \right)  \d p\d x\nonumber \\
   &&  - \frac12 \int_\Omega\int_\Omega \psi(p-x) f(p,t) \nabla_x \cdot \left(\grad\K[f](x,t) f(x,t) \right)  \d p \d x\nonumber \\
   &=& -\int_\Omega\int_\Omega w(x) \grad\varphi(p-x) \cdot \left( \grad\K[f](p,t) f(p,t) \right)  \d p\d x\nonumber \\
   && + \frac12 \int_\Omega\int_\Omega \grad\psi(p-x) f(x,t) \cdot \left(\grad\K[f](p,t) f(p,t) \right)  \d p\d x\nonumber \\
   &&  - \frac12 \int_\Omega\int_\Omega \grad\psi(p-x) f(p,t) \cdot \left(\grad\K[f](x,t) f(x,t) \right)  \d p \d x\nonumber \\
   &=& - \int_\Omega\int_\Omega \left( w(x) \grad\varphi(p-x) - \grad\psi(p-x) f(x,t) \right)
               \cdot \left(\grad\K[f](p,t) f(p,t) \right)  \d p\d x\nonumber \\
   &=& - \int_\Omega \left|\grad\K[f](p,t) f(p,t)\right|^2 f(p,t) \d p \leq 0 \,, \label{triangle}
\end{eqnarray}
where we have used the symmetry of $\psi$, i.e., $\grad\psi(p-x) = -\grad\psi(x-p)$,
and integration by parts, where  the boundary
terms vanish due to the compact support of $f$. 
(In the case $\Omega = \mathbb T^d$
due to the periodicity of the interaction kernels.)
\endproof

In the rest of this section, we study the question whether
the classical solution $f$ of~\eqref{kinetic}
tends for $t\to\infty$ to an equilibrium $f^*$,
characterized by the condition $\left|\grad\K[f^*](p) f^*(p)\right|^2 f^*(p) = 0$ a.e. on $\Omega$,
which stems from setting the right-hand side of~\eqref{triangle} equal to zero. 
Without loss of generality (see Remark~\ref{rem:3}), 
we restrict our attention to the case $\phi=\psi$;
then we have
\begin{equation} \label{phi=psi}
 \nabla {\cal K}[f](p,t) = \int_{\Omega} (w(x) - f(x,t)) \nabla \psi(p-x) \d x = \left( (w-f(\cdot,t))*\nabla \psi \right) (p).
\end{equation}
Moreover, we adopt the assumption that $f$ is a classical solution of~\eqref{kinetic}
in the sense of Lemma~\ref{lemma_1}.

First we show that the energy functional~\eqref{cont energy}
is bounded from below.

\begin{lemma} \label{lemma:boundedness}
Let $\Phi: [0,\infty) \rightarrow \R$ be concave and monotone increasing and $\Phi(0)$ be finite.
Then $\E[f]$ with $\phi=\psi=\Phi(| \cdot|)$ is bounded from below for all 
$f\in\mathcal M^1(\Omega)$.
\end{lemma}

\begproof
Without loss of generality, we can assume that $\Phi(0) = 0$.
Concavity and monotonicity of $\Phi$ imply its subadditivity~\cite{TSGSW},
which implies further
\[
    \Phi(|p-x|) \leq \Phi(|p|) + \Phi(|x|) \,,\quad
    \Phi(|p-x|) \geq \Phi(|p|) - \Phi(|x|) \,.
\]
Therefore, since $w, f \ge 0$ and $\int_\Omega w(p) \d p = \int_\Omega f(p) \d p = 1$, we obtain
\[
\E[f] &\geq& \int_\Omega \int_\Omega \left[ w(x) (\Phi(|p|)-\Phi(| x |)) - \frac12 f(x)(\Phi(|p|)+\Phi(|x|)) \right] f(p) \d p\d x \\
      &=& - \int_\Omega w(x) \Phi(|x|) \d x \,,
\]
which is the announced boundedness from below.
\endproof

\begin{remark}
Note that if $\Phi(r)$ grows as $r^\tau$ with $\tau > 2$ as $r \to\infty$
and $w$ is compactly supported,
then $\E$ is not bounded from below. This can be seen as follows:
Define $f_q(x) := \frac12 (\delta(x-q) + \delta(x+q))$ for $q\in\R^d$.
Then, due to the compact support of $w$, we get
\[
   \E[f_q] &=& \frac12 \int w(x) \left( \Phi(|q-x| + \Phi(|q+x|_2) \right) \d x - \frac14 \Phi(2|q|)\\
           &\sim& |q|^\tau - \frac14 (2 |q |)^\tau \; = \; (1 - \frac14 2^\tau)  |q|^\tau \to -\infty 
           \mbox{ as } |q| \to\infty \,.         
\]
\end{remark}

Now, since by Lemma \ref{lemma:classical_soln} the energy $\E[f(x,t)]$ is nonincreasing as a function of time
and by Lemma \ref{lemma:boundedness} bounded below, the limit of $\E[f(x,t)]$ as $t\to\infty$ exists and is finite.

Due to the boundedness of $f(\cdot,t)$ in the space of Radon measures,
there exists a sequence $t_j\to\infty$ 
and a Radon measure $f^*\in\mathcal M^1(\Omega)$
such that $f(\cdot,t_j) \to f^*$ weakly-* as $t_j\to\infty$.
Since, as assumed, $f$ is a classical solution with $f \in C^1(\Omega\times\mathbb R_+)$,
we also have $\E[f(\cdot,t_j)] \to \E[f^*]$ as $t_j\to\infty$,
and, consequently, $\E[f(\cdot,t)] \to \E[f^*]$ as $t_j\to\infty$.
By \eqref{triangle}, the equilibrium $f^*$ is characterized by the condition
\(  \label{equilibrium}
   \left|\grad\K[f^*](p) f^*(p)\right|^2 f^*(p) = 0 \quad \mbox{a.e. on } \Omega.
\)

By \eqref{phi=psi} the choice $f^*\equiv w$ is always a solution of~\eqref{equilibrium}.
This corresponds to the intuitive expectation that, as we let the number of particles $N$
tend to infinity, we should recover the profile $w$ in the long-time limit,
regardless of the initial distribution of particles.
However, the question whether the choice $f^*\equiv w$ is the unique solution of~\eqref{equilibrium}
in the class of probability measures seems to be rather nontrivial.
Although we believe that the affirmative indeed holds for a broad class
of interaction potentials, we are so far only able to provide a proof for the special case $\Omega=\R$
and $\phi(\cdot)=\psi(\cdot) = |\cdot|$.
Unfortunately, this case does not match our assumptions on the smoothness of $\phi$ and $\psi$
made in Theorem~\ref{theorem:mean-field} and Lemmata~\ref{lemma_1} and~\ref{lemma:classical_soln}.
However, the weak formulation~\eqref{weak} of equation~\eqref{kinetic} perfectly makes sense
if we insert the distributional derivative $\phi'(\cdot) = \psi'(\cdot) = \sign(\cdot)$ into~\eqref{K},
as long as $w, f\in L^1(\Omega)$, since then the integrals
\[
    \int_\Omega w(x) \phi'(p-x) \d x = \int_\Omega w(x) \sign(p-x) \d x 
\]
and
\[
    \int_\Omega f(x) \psi'(p-x) \d x = \int_\Omega f(x) \sign(p-x) \d x
\]
are well defined for all $p\in\R$ and uniformly bounded. Consequently, we can formulate the following Lemma:

\begin{lemma} \label{lemma:special_case}
Let $\Omega=\R$ and $\phi(\cdot)=\psi(\cdot) = |\cdot|$.
Let $w \ge 0$ be compactly supported in $\Omega$ and such that $\int_\Omega w(x) \d x = 1$.
Then the solution $f^* \equiv w$ of~\eqref{equilibrium} is unique in the class
$\mathcal{X} := \left\{ f\in L^1(\Omega): f\geq 0 \mbox{ a.e. on } \Omega, \int_\Omega f(x) \d x = 1 \right\}$.
\end{lemma}

\begproof
Let $f\in\mathcal{X}$ fulfill~\eqref{equilibrium}, which can be recast as
\begin{equation} \label{condition}
\left| \int_\Omega (w(x) - f(x)) \sign(p-x) \d x \right|^2 \, f(p) = 0 \quad \mbox{a.e. on } \Omega.
\end{equation}
Due to the normalization $\int_\Omega w(x) \d x = \int_\Omega f(x) \d x = 1$, we have
\begin{eqnarray*}
   G(p) &:=& \int_\Omega (w(x) - f(x)) \sign(p-x) \d x\\
        &=& \int_{-\infty}^p w(x) \d x - \int_p^{\infty} w(x) \d x - \int_{-\infty}^p f(x) \d x + \int_p^{\infty} f(x) \d x \\
        &=& 2W(p) - 2F(p) \,,
\end{eqnarray*}
where we denoted
\[
    W(p) := \int_{-\infty}^p w(x) \d x \qquad\mbox{and}\qquad F(p) := \int_{-\infty}^p f(x) \d x \,.
\]
The condition~\eqref{condition} can then be rewritten as $G(p)f(p)=0$ for almost all $p\in\Omega$;
let us note that $G$ is a continuous function.
By assumption, there exist real numbers $\alpha$ and $\omega$ such that $\supp w \subset [\alpha,\omega]$.
We prove that $G\equiv 0$ in three steps:
\begin{itemize}
\item First we show that $G(p) \equiv 0$ for all $p\leq\alpha$.
For a contradiction, let us assume that there exists a $p_0 \leq\alpha$ such that $G(p_0) \neq 0$.
Then we have $G(p_0) = 2W(p_0) - 2F(p_0) = -2F(p_0) < 0$ and, therefore,
there exists a set $S \subset (-\infty,p_0)$ of positive Lebesgue measure,
such that $f > 0$ almost everywhere on $S$.
But then $F > 0$ almost everywhere on $S$, and, consequently, $G = -2F < 0$ almost everywhere on $S$,
a contradiction to $G(p)f(p) = 0$.
\item Using a symmetry argument, we can prove that 
\begin{eqnarray*}
   G(p) = - 2\int_p^{\infty} w(x) \d x + 2\int_p^{\infty} f(x) \d x \equiv 0
\end{eqnarray*}
for all $p\geq\omega$.
\item Finally, we prove that $G(p)\equiv 0$ also for all $p\in (\alpha,\omega)$.
By the continuity of $G$, 
the set 
${\cal S} := \{p \in (\alpha,\omega): G(p) \neq 0\}$ is open.
Since~\eqref{condition} dictates that $f(p)=0$ almost everywhere on ${\cal S}$,
we have for every $p_0\in {\cal S}$ and $\delta>0$ small enough,
\[
    G(p_0+\delta) = G(p_0) + 2 \int_{p_0}^{p_0+\delta} w(x) \d x \geq G(p_0) \,.
\]
Therefore, $G$ is nondecreasing on ${\cal S}$, and, thus, $G$ is nondecreasing everywhere on $(\alpha,\omega)$.
Since $G(\alpha) = G(\omega) = 0$ and $G$ is continuous, we conclude that $G\equiv 0$.
\end{itemize}
We finish the proof by observing that
\[
    \int_{-\infty}^p (w(x) - f(x)) \d x = 0 \qquad\mbox{for all } p\in\R
\]
implies $f=w$ almost everywhere on $\R$.
\endproof

\section{Numerical Examples} \label{sec:numerics}

\begin{figure}[htb]
{
\begin{tabular}[h]{ll}
\resizebox*{0.48\linewidth}{!}{\includegraphics{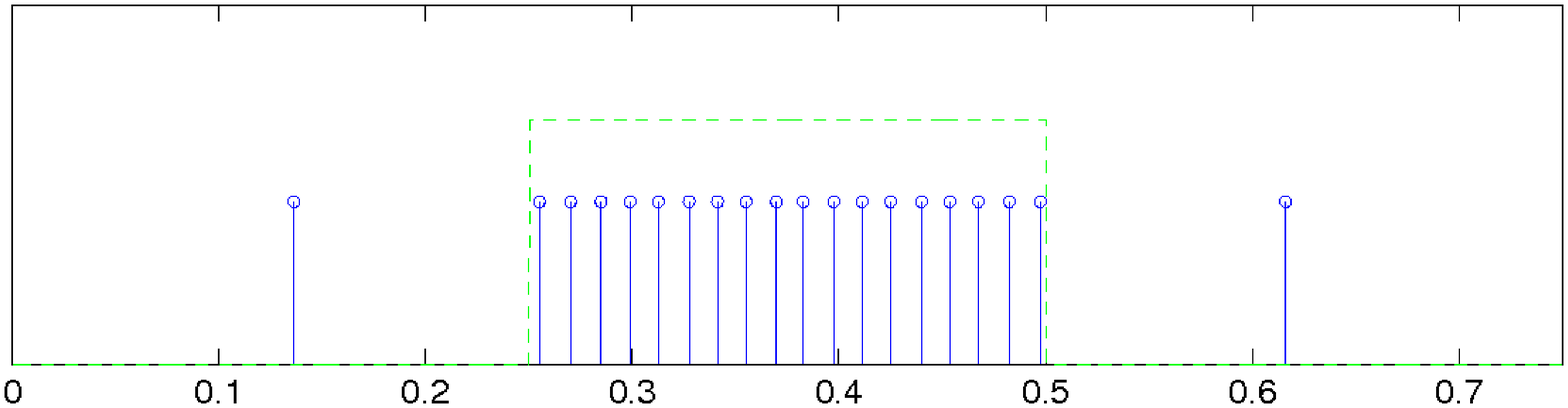}} &
\resizebox*{0.48\linewidth}{!}{\includegraphics{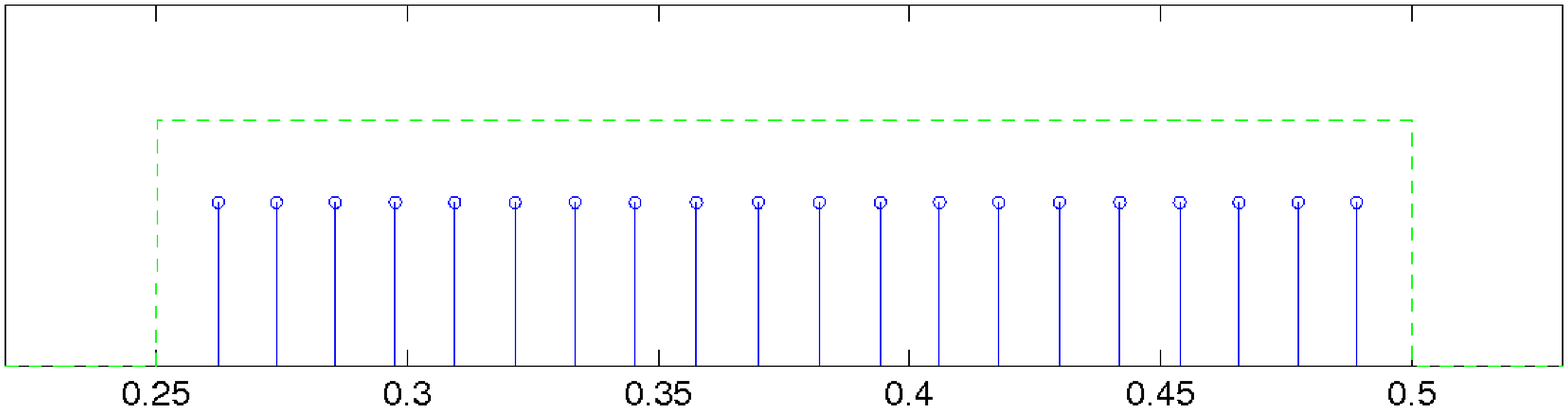}} \\
\resizebox*{0.48\linewidth}{!}{\includegraphics{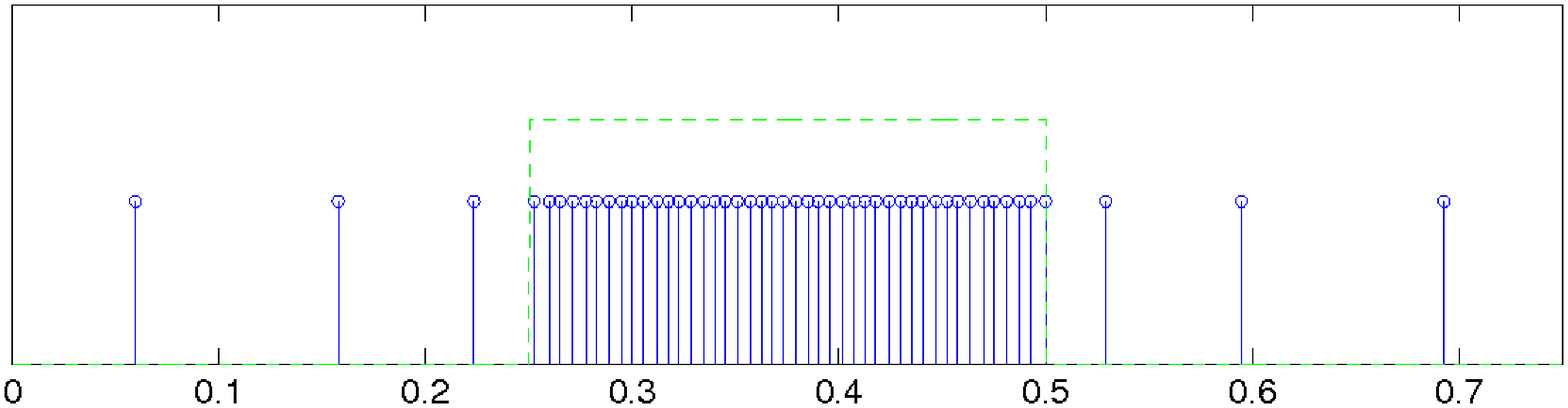}} &
\resizebox*{0.48\linewidth}{!}{\includegraphics{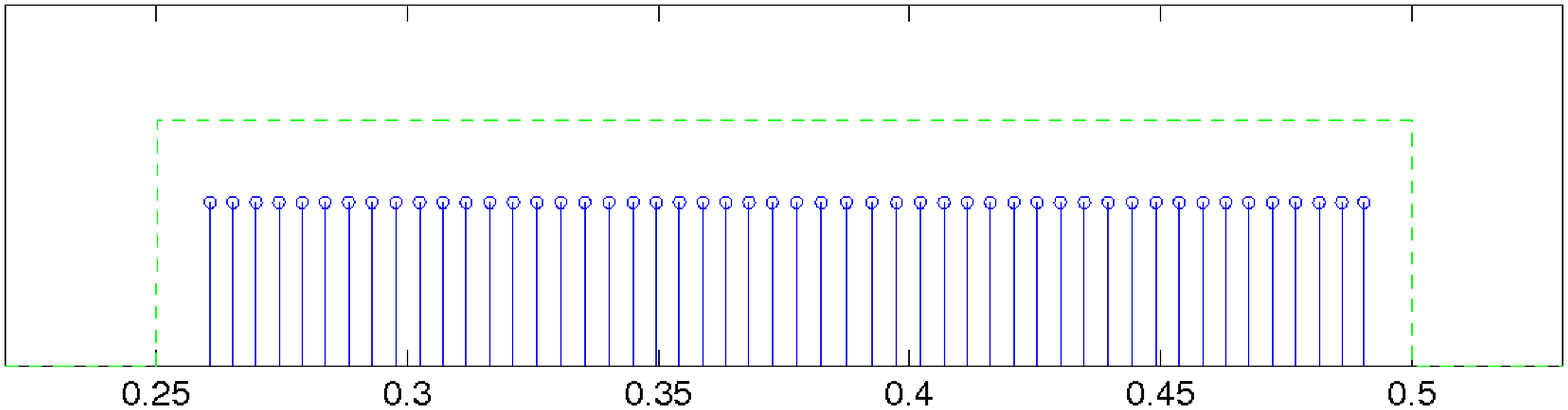}} \\
\resizebox*{0.48\linewidth}{!}{\includegraphics{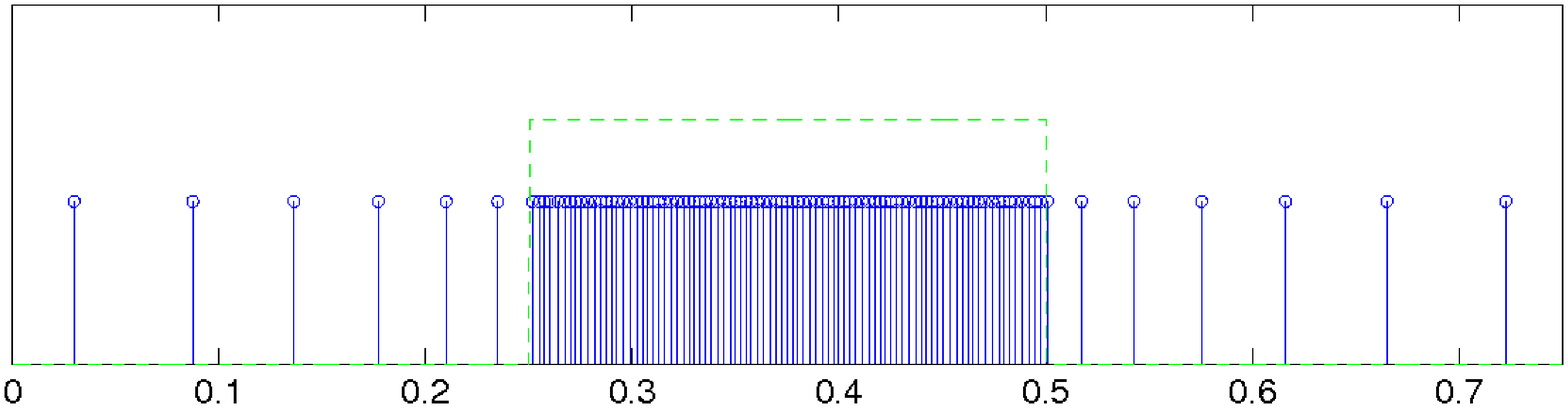}} &
\resizebox*{0.48\linewidth}{!}{\includegraphics{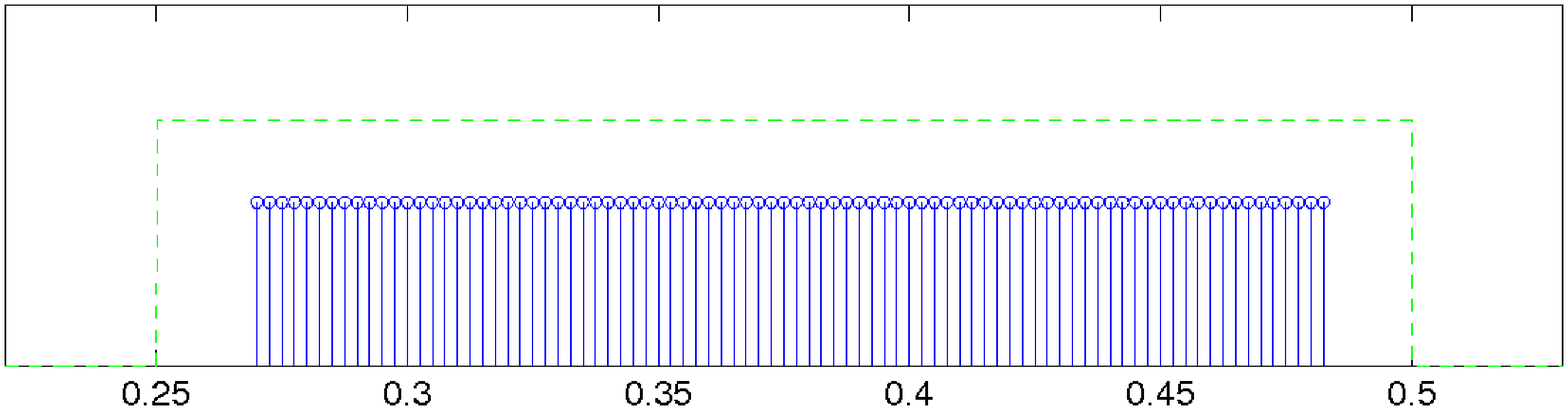}}
\end{tabular}\par}
\caption{Steady states of the discrete particle system~\eqref{ODE}
with $N=20, 50, 100$ particles (first, second and, resp., third row)
in the smoothing (left panels) and sharpening case (right panels).
The datum $w$ is visualized with the dashed line.
Note the different horizontal axis limits in the smoothing and sharpening case.
}\label{fig:num_discrete}
\end{figure}

In this section we present few numerical results
for the discrete particle system~\eqref{kernel_psi_phi}
and the mean-field limit~\eqref{kinetic}.
We consider two cases:
\begin{itemize}
\item Smoothing case: $\phi(s) = |s|^{1.1}$ and $\psi(s) = |s|$
\item Sharpening case: $\phi(s) = |s|$ and $\psi(s) = |s|^{1.1}$
\end{itemize}
In both cases, we consider the 1D full-space setting $\Omega = \R$
with the datum $w = 4\chi_{[0.25,5]}$.
For the discrete particle system, we use $N=20, 50$ and $100$ particles.
We integrate the  ODE system~\eqref{ODE} in time using
the explicit Euler method until the steady state (which does not depend on the initial condition).
The results for the smoothing and sharpening cases are shown in Fig.~\ref{fig:num_discrete}.

\noindent
\vspace{0.3cm}

\begin{figure}[htb]
{\centering
\begin{tabular}[h]{cc}
\resizebox*{0.45\linewidth}{!}{\includegraphics{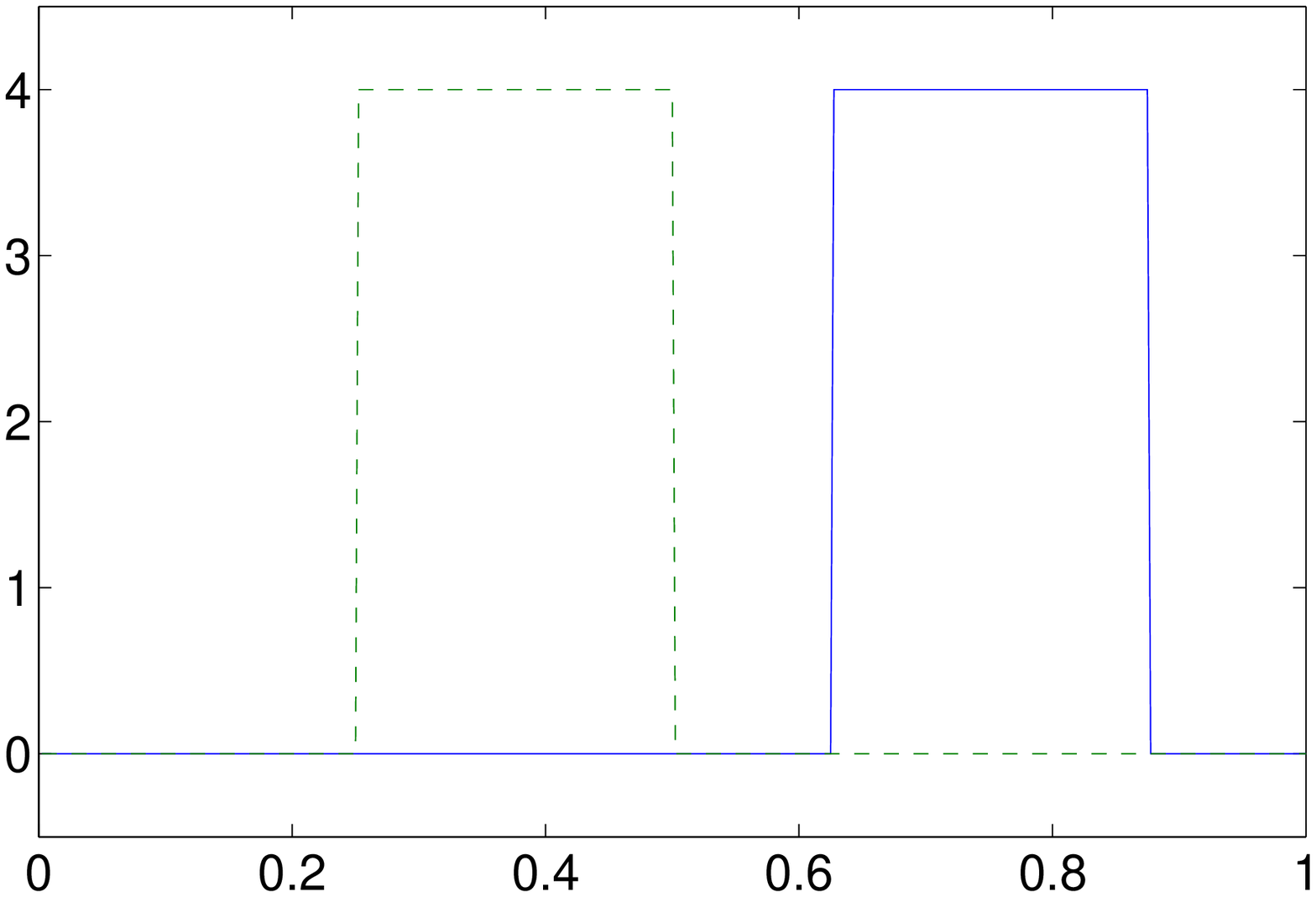}} &
\resizebox*{0.45\linewidth}{!}{\includegraphics{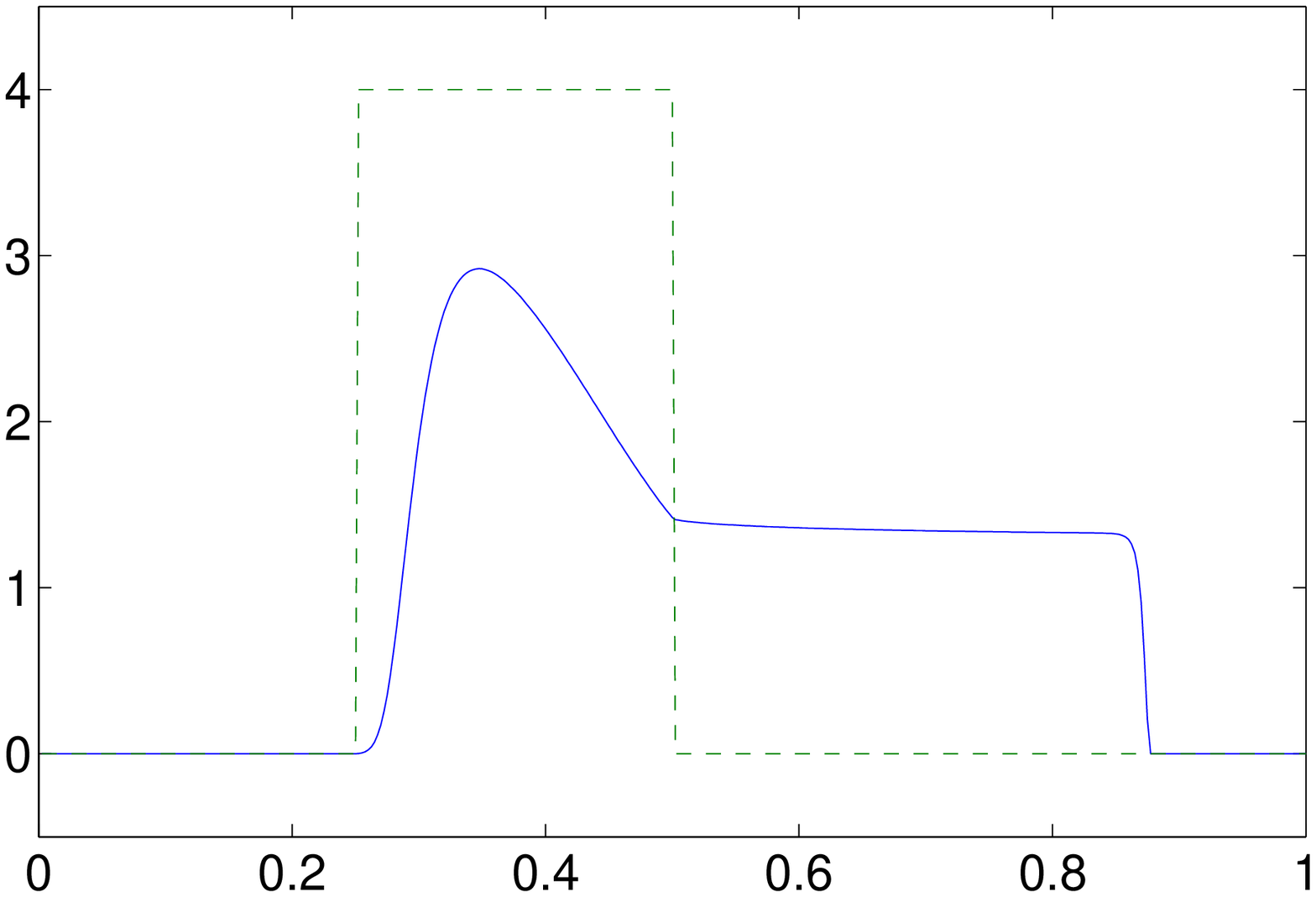}} \\
$t=0.0$ & $t=1.5$ \\
\resizebox*{0.45\linewidth}{!}{\includegraphics{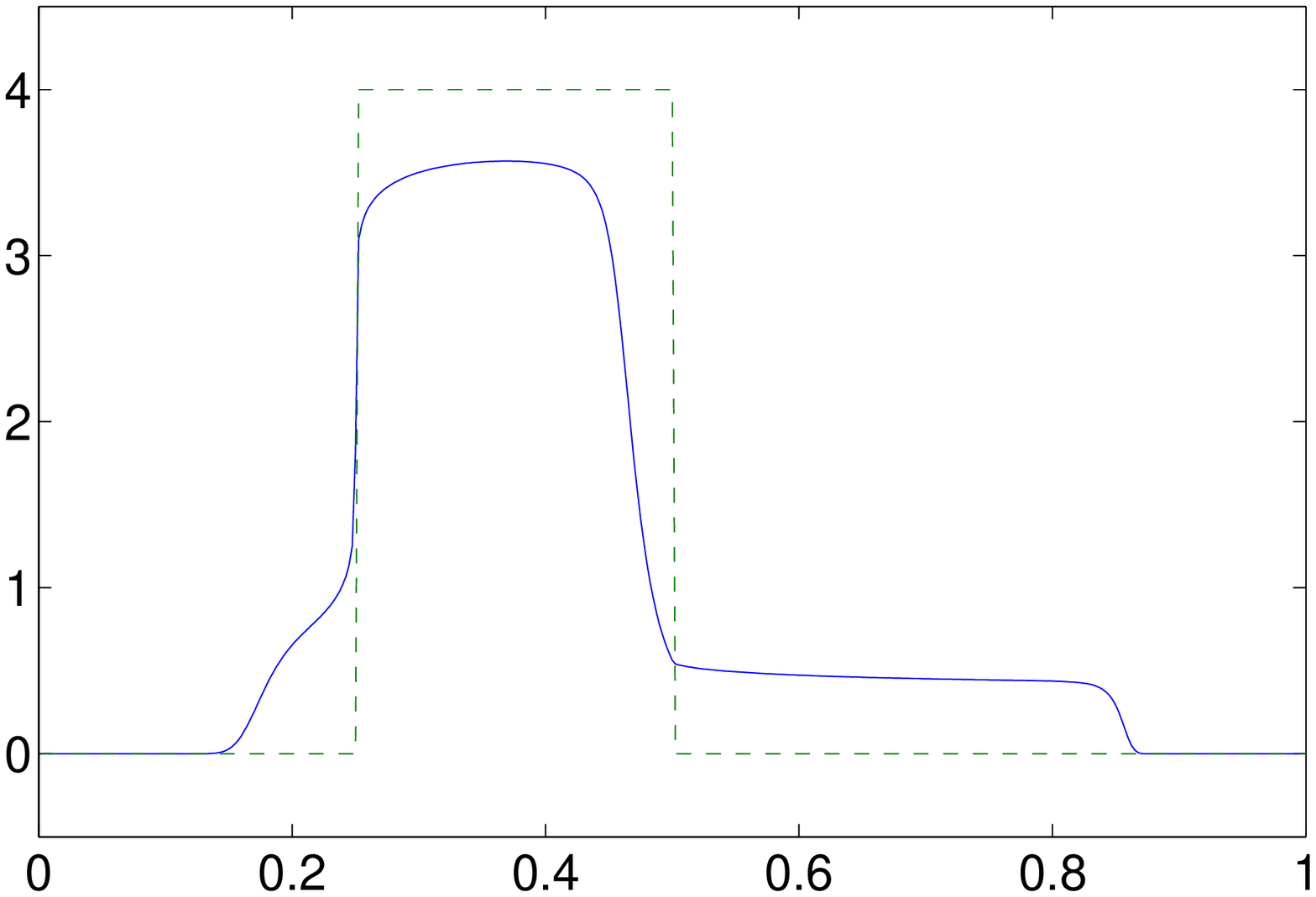}} &
\resizebox*{0.45\linewidth}{!}{\includegraphics{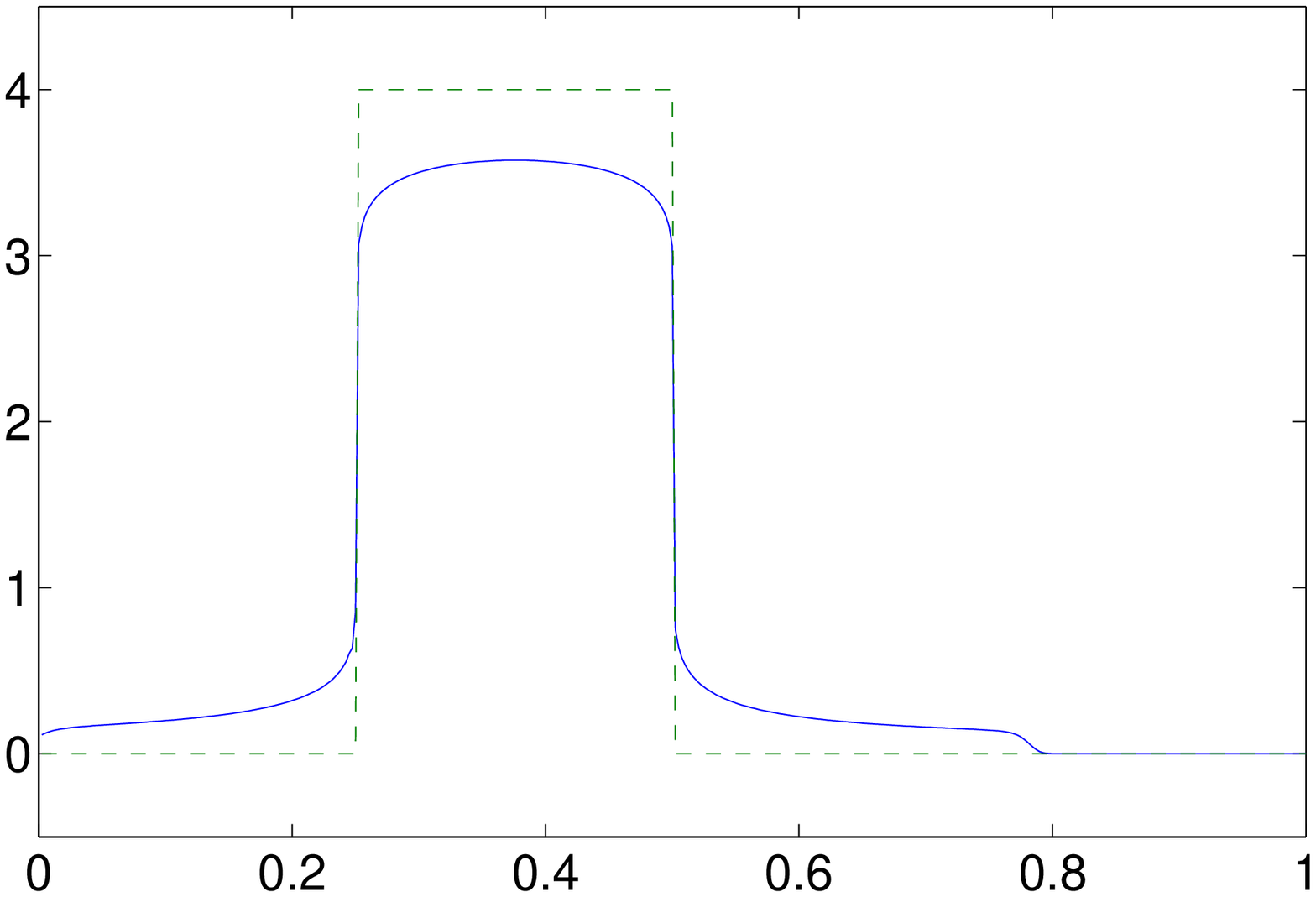}} \\
$t=3.0$ & $t=5.0$
\end{tabular}\par}
\caption{Smoothing case of the mean-field limit~\eqref{kinetic}.
The solid line represents the solution $f$, the dashed line the data $w$.
The upper left panel shows the initial condition, the lower right panel is the steady state.
}\label{fig:num3}
\end{figure}

\begin{figure}[htb]
{\centering \begin{tabular}[h]{cc}
\resizebox*{0.45\linewidth}{!}{\includegraphics{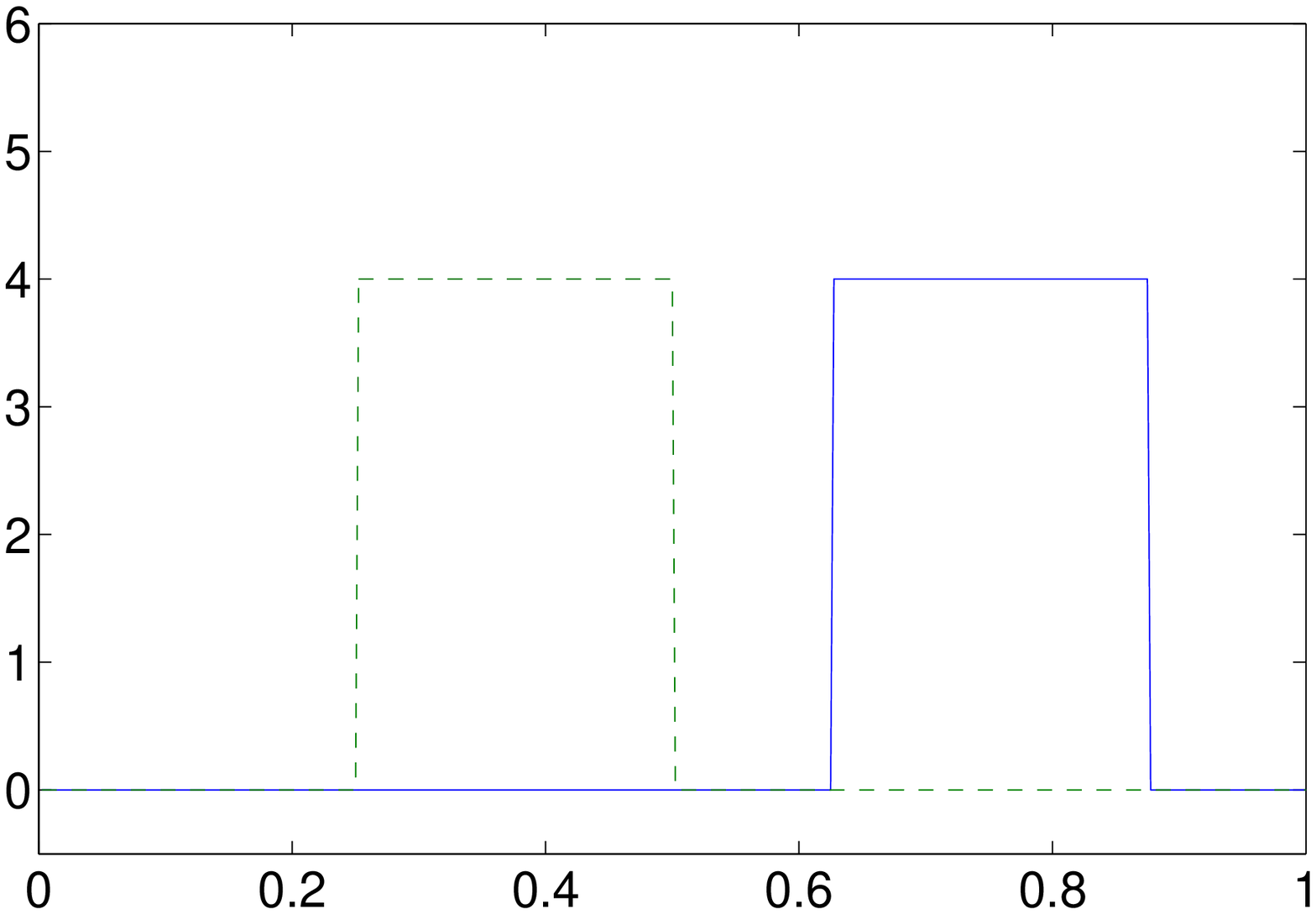}} &
\resizebox*{0.45\linewidth}{!}{\includegraphics{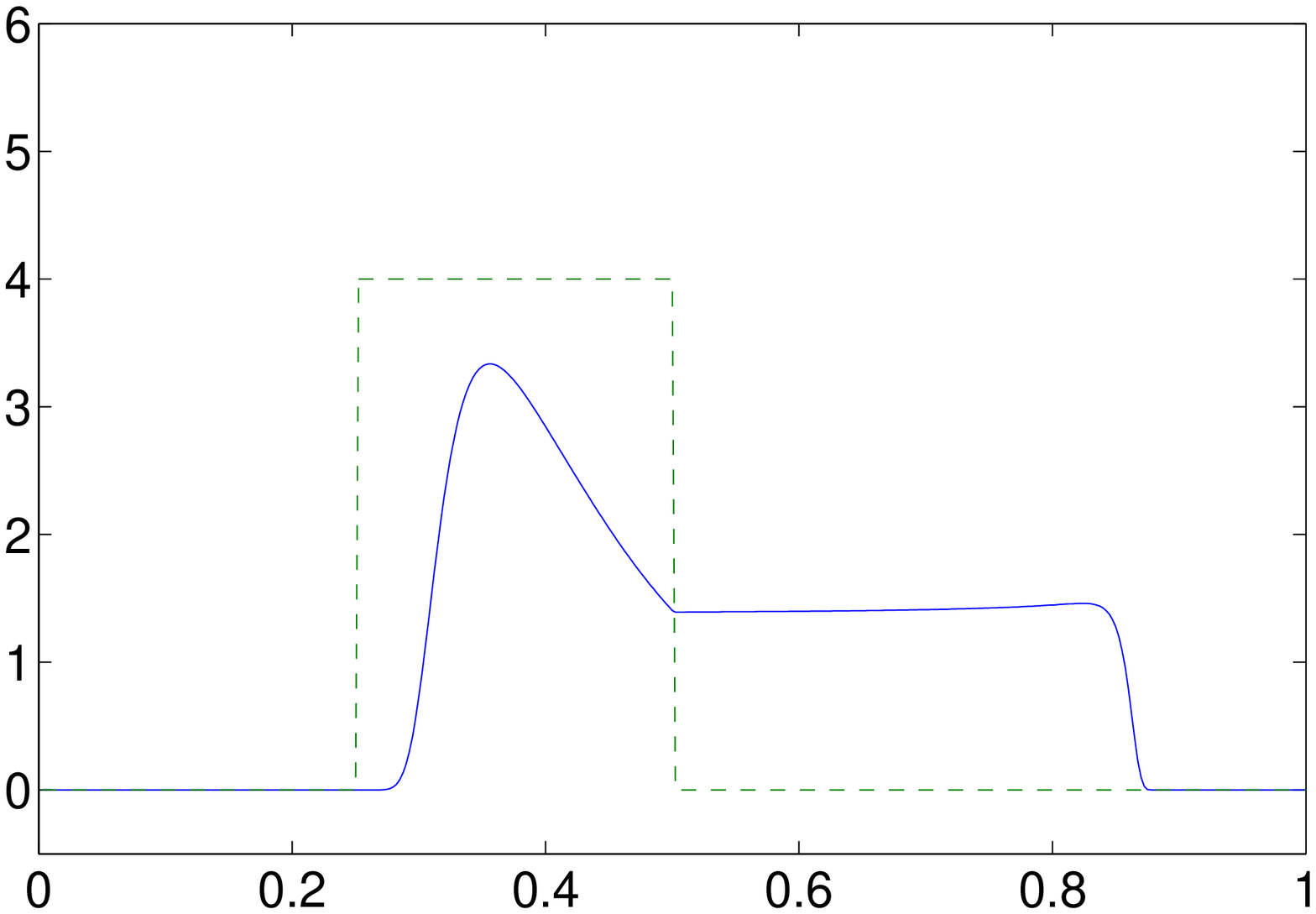}} \\
$t=0.0$ & $t=1.5$ \\
\resizebox*{0.45\linewidth}{!}{\includegraphics{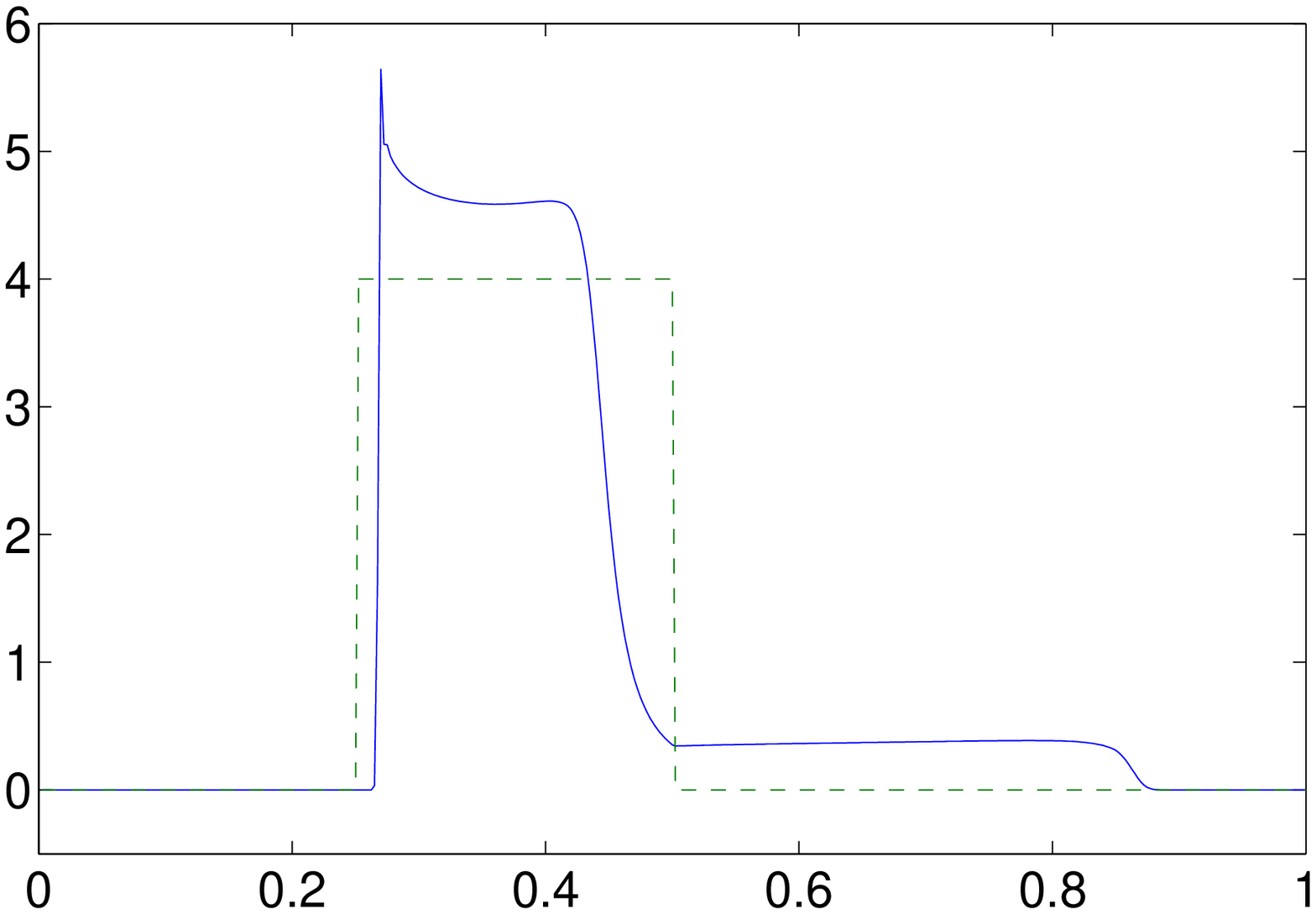}} &
\resizebox*{0.45\linewidth}{!}{\includegraphics{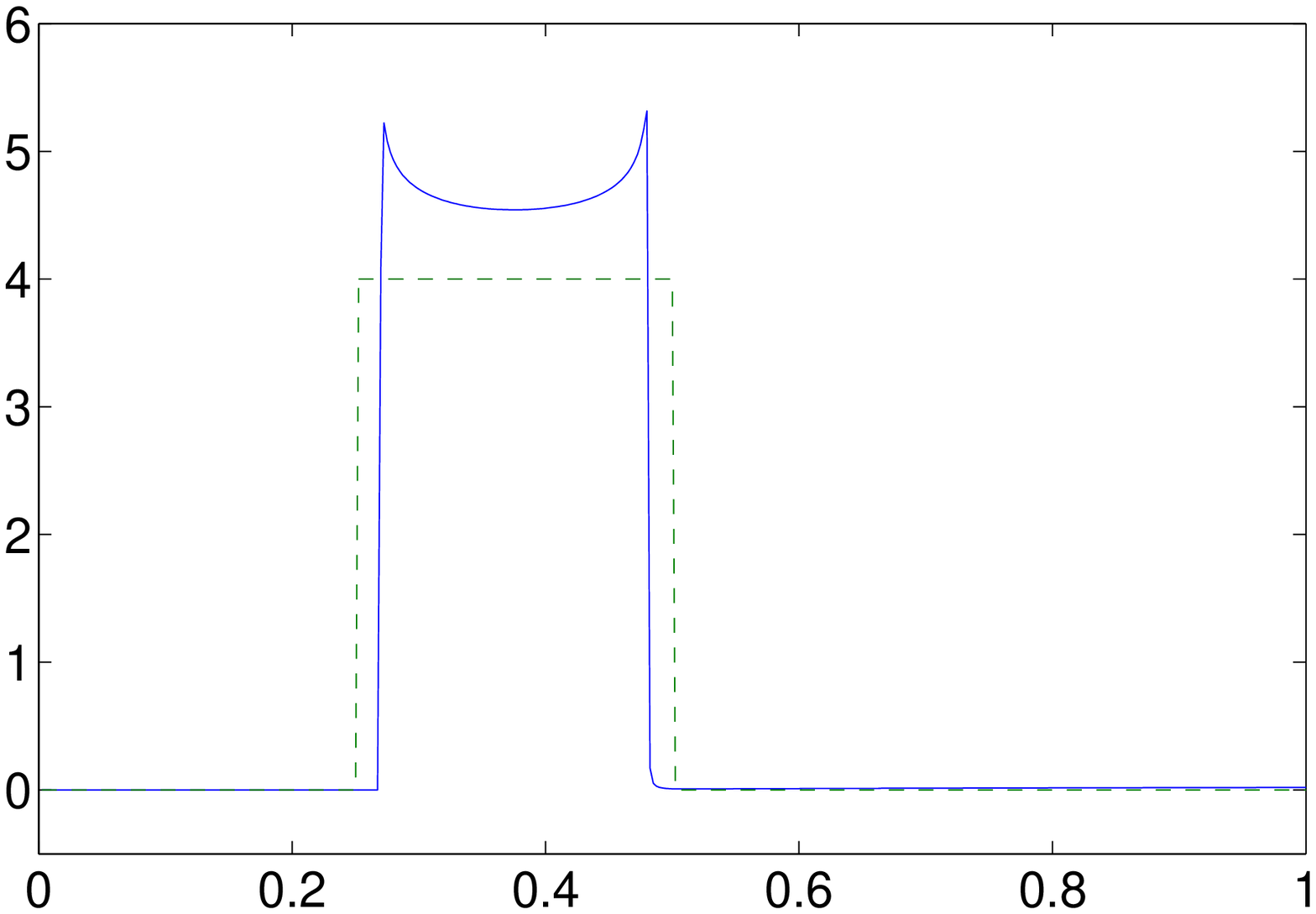}} \\
$t=3.5$ & $t=9.5$
\end{tabular}\par}
\caption{Sharpening case of the mean-field limit~\eqref{kinetic}.
The solid line represents the solution $f$, the dashed line the data $w$.
The upper left panel shows the initial condition, the lower right panel is the steady state.
}\label{fig:num4}
\end{figure}

For the mean-field limit~\eqref{kinetic}, we impose the initial condition
$f_0 = 4\chi_{[0.65,0.9]}$.
We discretize~\eqref{kinetic} using the semi-implicit finite difference method
with upwinding in space and explicit Euler method in time.
Snapshots of the solutions
are shown in Fig.~\ref{fig:num3} for the smoothing case
and in Fig.~\ref{fig:num4} for the sharpening case.

As we expect, the solutions converge to some steady states as $t\to\infty$
in both the discrete and mean-field cases.
In the smoothing case, the equilibrium profile is a smoothed version of the data $w$,
while in the sharpening case the equilibrium is an anti-smoothed version of $w$.
It is interesting to compare these results with the discrete particle calculation
from the previous example. Indeed, the steady state particle distributions
shown on the left panels of Fig.~\ref{fig:num_discrete} can be regarded as approximations
of the steady state on Fig.~\ref{fig:num3}, and the same hold for the right panels of
Fig.~\ref{fig:num_discrete} and the steady state of Fig.~\ref{fig:num4}.

\noindent
\vspace{0.3cm}

\noindent

\begin{figure}[htb]
{\centering \begin{tabular}[h]{cc}
\resizebox*{0.45\linewidth}{!}{\includegraphics{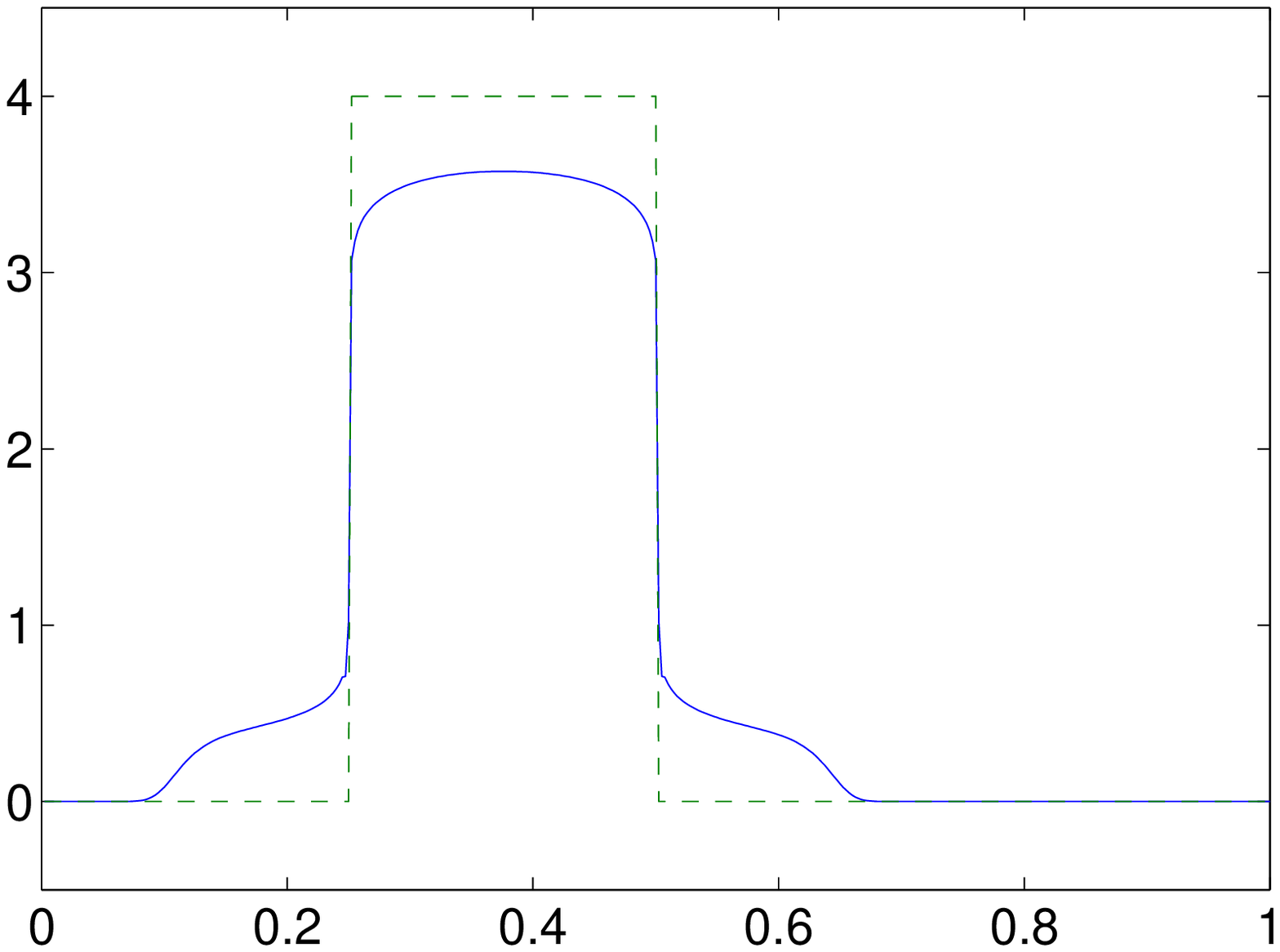}} &
\resizebox*{0.45\linewidth}{!}{\includegraphics{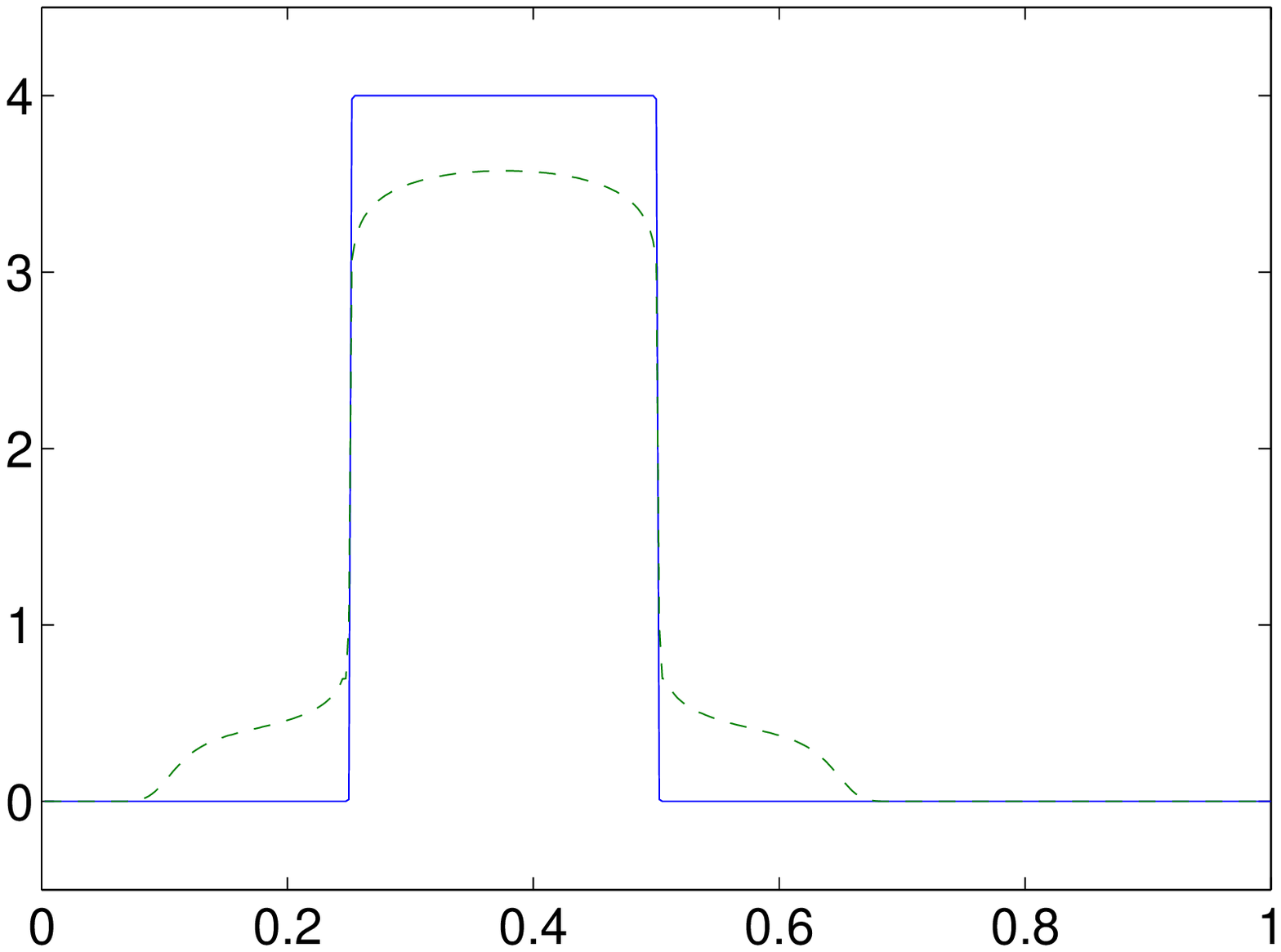}}
\end{tabular}\par}
\caption{Left panel: The original image $w$ (dashed line) and its blurred version $\tilde w=f(t=2)$ (solid line).
Right panel: The result of deblurring of $\tilde w$ (dashed line), given
by the steady state $f(\cdot,t=30)$ of~\fref{kinetic}, solid line.
}\label{fig:num5}
\end{figure}

In our last example, we show how the kinetic equation~\eqref{kinetic} can be used
for deblurring of images in the case when the blurring kernel is known.
As our ``image'' we take again the datum $w = 4\chi_{[0.25,5]}$,
and blur it by applying $200$ time-steps of~\eqref{kinetic}
with $\phi(s) = |s|^{1.1}$ and $\psi(s) = |s|$ (smoothing case),
the time-step length is $0.01$.
For simplicity, we use the initial datum $f_0:=w$.
The result of this blurring process, $g= f(\cdot,t=2)$, is shown on Fig.~\ref{fig:num5}, left panel.
Then, we perform deblurring of the ``image''
by applying~\eqref{kinetic} with reversed interaction potentials,
i.e., $\phi(s) = |s|$ and $\psi(s) = |s|^{1.1}$,
and with $w:=g$.
Again, for simplicity, we take $f_0:=w=g$ as the initial condition for $f$.
We let~\eqref{kinetic} evolve until steady state, which is reached around $t=30$,
and plot it in the right panel of Fig.~\ref{fig:num5}.
We see that $f(\cdot,t=30)$ is indistinguishable from the original, un-blurred image $w$;
the relative difference between them in the $L^1$-norm is approximately $0.6\%$.

\noindent
\vspace{0.3cm}

Let us mention that numerical implementation of~\eqref{kinetic} in the spatially 2D setting,
which potentially might be of interest for application in image processing (deblurring of images),
is quite a demanding task. The main reason is the high numerical cost,
caused by the necessity of evaluation of the convolution $f*\grad\psi$ in each time step,
which in general takes $\mathcal O(N^4)$ multiplications if the grid consists of $N^2$ points.
A possible speed-up of this operation can be achieved by using fast multipole expansion or FFT-based methods,
see for instance~\cite{Fenn}. We postpone this task for future work.

\subsubsection*{Acknowledgments}
Massimo Fornasier and Jan Ha\v{s}kovec  acknowledge
 the financial support provided by the START award ``Sparse Approximation and Optimization in High Dimensions'' no. FWF~Y~432-N15
of the Fonds zur F\"orderung der wissenschaftlichen Forschung (Austrian Science Foundation).


\end{document}